\documentclass[a4paper,12pt]{article}
\usepackage{amsmath,amsthm,amsfonts,latexsym,amscd}
\newtheorem{theorem}{Theorem}[section]
\newtheorem{lemma}[theorem]{Lemma}
\newtheorem{corollary}[theorem]{Corollary}
\newtheorem{proposition}[theorem]{Proposition}

\theoremstyle{definition}
\newtheorem{definition}[theorem]{Definition}
\newtheorem{example}[theorem]{Example}
\newtheorem{notation}[theorem]{Notation}

\newtheorem{remark}[theorem]{Remark}

\newcommand{\R}{\mathbb{R}}
\newcommand{\C}{\mathbb{C}}
\newcommand{\N}{\mathbb{N}}
\newcommand{\Z}{\mathbb{Z}}

\newcommand{\HH}{\mathbb{H}}

\DeclareMathOperator{\id}{id}
\newcommand{\boundary}[1]{\partial#1}

\newcommand{\abs}[1]{\lvert#1\rvert} 
\newcommand{\norm}[1]{\lVert#1\rVert}

\DeclareMathOperator{\supp}{supp}   
\DeclareMathOperator{\clos}{clos}  
\DeclareMathOperator{\im}{im}      
\DeclareMathOperator{\vol}{vol}    

\DeclareMathOperator{\tr}{Tr}
\DeclareMathOperator{\pr}{pr}
\DeclareMathOperator{\Spur}{tr}    

{\catcode`@=11\global\let\c@equation=\c@theorem}

\DeclareMathOperator{\ordanTor}{T_{\text{\rm an}}}      
\DeclareMathOperator{\anTor}{T^{(2)}_{\text{\rm an}}}   
\DeclareMathOperator{\topTor}{T^{(2)}_{\text{\rm top}}} 
\newcommand{\Neumann}[1]{\mathcal{#1}}
\DeclareMathOperator{\Ex}{Ex}

\begin{document}


\setcounter{section}{0}
\typeout{--------------------- introduction -----------------------}

\title{$L^2$-torsion of hyperbolic manifolds of finite volume}
\author{ Wolfgang L\"uck \and Thomas Schick}
\maketitle

\begin{abstract}
Suppose $\overline{M}$ is a compact connected
odd-dimensional manifold with boundary, whose
interior $M$ comes with a complete
hyperbolic metric of finite volume. We will show
that the $L^2$-topological torsion of
$\overline{M}$ and the $L^2$-analytic torsion
of the Riemannian manifold $M$ are equal.
In particular, the $L^2$-topological torsion of
$\overline{M}$ is proportional to the
hyperbolic volume of $M$, with a constant of
proportionality which depends only on the
dimension and which is known to be nonzero in
odd dimensions \cite{Hess(1998)}. In dimension 3 this proves
the conjecture
\cite[Conjecture 2.3]{Lueck (1994a)} or
\cite[conjecture 7.7]{Lott-Lueck (1995)} which gives a
complete calculation of the $L^2$-topological torsion of
compact $L^2$-acyclic $3$-manifolds which admit
a geometric JSJT-decomposition.

In an appendix we give a counterexample to an extension of the Cheeger-M\"uller
theorem to manifolds with boundary: if the metric is not a product
near the boundary, in general analytic and topological torsion are not
equal, even if the Euler characteristic of the boundary vanishes.\\[2mm]
AMS-classification number: 58G11
\end{abstract}

\setcounter{section}{-1}
\section{Introduction}\label{Introduction}

In this paper we study $L^2$-analytic
torsion of a compact connected manifold $\overline{M}$ with
boundary such that the interior $M$ comes with
a complete hyperbolic metric of finite
volume.

\begin{notation} \label{notation for the ends}
Let $m$ be the dimension of $M$.
From hyperbolic geometry we know
\cite[Chapter  D3]{Benedetti-Petronio (1992)}
that $M$ can be written as
\[ M = M_0\cup_{\boundary M_0} E_0, \]
where $M_0$ is a compact manifold with
boundary and $E_0$ is a finite disjoint
union of hyperbolic ends $[0,\infty)\times F_j$.
Here each $F_j$ is a closed
flat manifold and the metric on the end is the warped product
 \[ du^2+e^{-2u}dx^2 \]
with $dx^2$ the metric of $F_j$.
Of course, we can make the ends smaller and
also write
\[ M= M_R\cup_{\boundary M_R} E_R \qquad R\ge 0,\]
where $E_R$ is the subset of $E_0$ consisting of the components
$[R,\infty)\times F_j$. We define
\[ T_R=M_{R+1}\cap E_R\qquad R\ge 0.\]
Denote by $\widetilde{M}_R$, $\widetilde{E_R}$, $\ldots$
the inverse images of $M_R$, $E_R$,
$\ldots$ under the universal covering
map $\widetilde{M} \longrightarrow M$.
Correspondingly, for the universal covering we have
\[ \widetilde{M} =\widetilde{M}_R\cup_{\boundary\widetilde{M}_R}
\widetilde{E_R}. \]
For $\widetilde{E_R}$ we get and fix
coordinates such that each component is
\[ [R,\infty)\times \R^{m-1}\quad\text{with warped product metric
$du^2+ e^{-2u}dx^2$}, \]
where $dx^2$ is the Euclidean metric on $\R^n$.\qed
\end{notation}

Each $M_R$ is a compact connected Riemannian
manifold with boundary which is  $L^2$-acyclic (see Corollary
\ref{L2_acyclic}). Its absolute
$L^2$-analytic torsion $\anTor(M_R)$ is defined.
The manifold $M$ has no boundary and finite volume
and its $L^2$-analytic torsion is defined although
it is not compact (Remark \ref{definitions make also sense in the
hyperbolic finite volume case}). We will recall the notion of
$L^2$-analytic torsion in Section
\ref{Review of L^2-analytic torsion and strategy of proof}.
The main result of this paper is

\begin{theorem}\label{Tor_converg}
Let $\overline{M}$ be a compact connected manifold with
boundary such that the interior $M$ comes with
a complete hyperbolic metric of finite wolume. Suppose that
the dimension $m$ of $M$ is odd. Then we get, if
$R$ tends to infinity
\[ \lim_{R \to \infty} \anTor(M_R) ~ = ~\anTor(M). \qed \]
\end{theorem}

\begin{remark}
 For compact Riemannian manifolds with boundary and with a product
 metric near the boundary,
 analytic torsion and topological torsion differ by $(\ln 2)/2$ times
 the Euler characteristic of the boundary. This is a result of L\"uck
 in the classical situation \cite{Lueck(1993)} and of Burghelea et.al
 \cite{Burghelea-Friedlander-Kappeler (1996a)} for the
 $L^2$-version. In particular, analytic torsion does not depend on the
 metric, as long as it is a product near the boundary and the manifold
 is acyclic or $L^2$-acyclic, respectively.\\
 In
 appendix \ref{metric_anomaly_example} we give examples which show that
 this is not longer the case for arbitrary metrics. This answers an
 old question of Cheeger's \cite[p.~281]{Cheeger (1979)}. Of course it also
 implies that the above extension of the Cheeger-M\"uller theorem to
 manifolds with boundary is not true in general.

 This requires additional care in the chopping and exhausting process
 described above.
\end{remark}

We will explain the strategy of proof for \ref{Tor_converg}
in Section
\ref{Review of L^2-analytic torsion and strategy of proof}. Next we discuss
consequences of Theorem \ref{Tor_converg} and put it into context
with known results.\par

We will explain in Section
\ref{Review of L^2-analytic torsion and strategy
of proof} that
the comparison theorem for $L^2$-analytic
and -topological torsion for
manifolds with and without boundary of Burghelea, Friedlander,
Kappeler and McDonald
\cite{Burghelea-Friedlander-Kappeler  (1996a),%
Burghelea-Friedlander-Kappeler McDonald (1996a)}
now implies

\begin{theorem}\label{TopTor=AnTor}
Let $\overline{M}$ be a compact connected manifold with
boundary such that the interior $M$ comes with
a complete hyperbolic metric of finite wolume.
Then
\[ \topTor(\overline{M})=\anTor(M). \qed\]
\end{theorem}

The computations of Lott \cite[Proposition 16]{Lott (1992a)}
and Mathai \cite[Corollary 6.7]{Mathai (1992)} for closed
hyperbolic manifolds extend directly to hyperbolic manifolds
without boundary and with finite volume since $\HH^m$ is homogeneous.
(Notice that we use for the analytic
torsion the convention in Lott which is twice the logarithm of the one of Mathai).
Hence Theorem \ref{TopTor=AnTor} implies 

\begin{theorem}\label{AnTor = Vol}
Let $\overline{M}$ be a compact connected manifold with
boundary such that the interior $M$ comes with
a complete hyperbolic metric of finite wolume.
Then there is a dimension constant $C_m$ such that
\[\anTor(M) ~ = ~ C_m \cdot \vol(M). \]
Moreover, $C_m$ is zero, if $m$ is even, and
$C_3 = \frac{-1}{3\pi}$
and $(-1)^m C_{2m+1}>0$. \qed
\end{theorem}
\begin{remark}
  The last statement is a result of Hess \cite{Hess(1998)}
  and answers the question of Lott
  \cite{Lott (1992a)} wether $C_{2m+1}$ is always nonzero.

  For closed manifolds, the proportionality follows directly from the
  computations of Lott and Mathai and the
  comparison theorem for $L^2$-analytical and topological torsion of
  Burghelea
  et.al. \cite{Burghelea-Friedlander-Kappeler McDonald (1996a)} (as is
  also observed there).
\end{remark}

Now Theorem \ref{TopTor=AnTor} and Theorem \ref{AnTor = Vol}
together with \cite[Theorem 2.1]{Lueck (1994a)} imply

\begin{theorem}\label{anTor for 3-manifolds}
Let $N$ be a compact connected orientable irreducible $3$-manifold with
infinite fundamental group possessing a
geometric JSJT-decomposition such that
the boundary of $N$ is empty or a
disjoint union of incompressible tori.
Then all $L^2$-Betti numbers of $N$ vanish and
all Novikov-Shubin invariants are positive. Moreover, we get
\[ \topTor(N) ~ = ~\frac{-1}{3\pi}
\cdot \sum_{i=1}^r \vol(N_i), \]
where $N_1$, \ldots $N_r$ are the hyperbolic pieces of
finite volume in the JSJT-decomposition of $N$.
In particular $\topTor(N)$ is $0$ if and only if
$N$ is a graph-manifold, i.e. there are no
hyperbolic pieces of finite volume in the JSJT-decomposition. \qed
\end{theorem}

Theorem \ref{anTor for 3-manifolds} has been conjectured
in \cite[Conjecture 2.3]{Lueck (1994a)} and
\cite[Conjecture 7.7]{Lott-Lueck (1995)}, where also the relevant notions
are explained. The JSJT-decomposition of an irreducible connected
orientable compact $3$-manifold $N$ with infinite fundamental group
is the decomposition of Jaco-Shalen and Johannson by a
minimal family of pairwise non-isotopic incompressible
not boundary-parallel embedded $2$-tori into Seifert pieces
and atoroidal pieces. If these atoroidal pieces are hyperbolic, then
the JSJT-decomposition is called geometric.
Thurston's Geometrization Conjecture
says that these atoroidal pieces are always hyperbolic.
This conjecture is known to be true if
 $N$ is Haken, for instance if $N$ has boundary or
its first Betti number is positive.

In this context we mention the combinatorial formula
for the topological torsion of a $3$-manifold $N$ as in Theorem 
\ref{anTor for 3-manifolds} which computes
this torsion just from a presentation of its fundamental
group without using any information about $N$
itself \cite[Theorem 2.4]{Lueck (1994a)}. Theorem
\ref{anTor for 3-manifolds} also implies that for such $3$-manifolds
 the $L^2$-torsion and the simplicial
volume of Gromov agree up to a non-zero
multiplicative constant \cite[section 2]{Lueck (1994a)}.

The paper is organized as follows:\\[4mm]
\begin{tabular}{ll}
\ref{Introduction}. & Introduction
\\
\ref{Review of L^2-analytic torsion and strategy of proof}.
& Review of $L^2$-analytic torsion and strategy of proof
\\
\ref{Analysis of the heat kernel}.
& Analysis of the heat kernel
\\
\ref{The large t summand of the torsion
corresponds to small eigenvalues}.
& The large t summand of the torsion
corresponds to small eigenvalues
\\
\ref{Spectral density functions}.
& Spectral density functions
\\
\ref{Sobolev- and L^2-complexes}.
& Sobolev- and $L^2$-complexes
\\
\ref{Spectral density functions for M_R}.
& Spectral density functions for $M_R$
\\
\ref{L^2-analytic torsion and variation of the metric}
& $L^2$-analytic torsion and variation of the metric
\\
\ref{metric_anomaly_example} & Examples for nontrivial metric anomaly\\
& References
\end{tabular}



\typeout{--------------------- section 1-----------------------}
\section{Review of $L^2$-analytic torsion and strategy of proof}
\label{Review of L^2-analytic torsion and strategy of proof}

In this section we recall the definition of $L^2$-analytic
torsion (compare 
\cite{Novikov-Shubin(1986a),Novikov-Shubin(1986b),Lott (1992a),Mathai (1992)}),
discuss the strategy and set the
stage for the proof of Theorem \ref{Tor_converg}, and prove Theorem
\ref{TopTor=AnTor}.

\begin{definition} \label{definition of normalized trace}
Let $N$ be a connected compact $m$-dimensional Riemannian manifold.
Let $\Delta_p[\widetilde{N}]$ be the Laplacian on
$p$-forms on the universal covering $\widetilde{N}$
considered as an unbounded
self-adjoint operator on $L^2$ (with absolute
boundary conditions, if the boundary is non-empty). Then
$e^{-t\Delta_p[\widetilde{N}]}$ is defined for
every $t>0$ and has a smooth kernel $e^{-t\Delta_p[\widetilde{N}]}(x,y)$.
(This kernel is invariant under the diagonal
$\Gamma :=\pi_1(N)$-operation on
$\widetilde{N}\times\widetilde{N}$
given by deck transformations). Let ${\cal F}$ be a fundamental domain
for the $\Gamma$-action. Now define the
{\em normalized trace}
\[ \tr_{\Gamma}e^{-t\Delta_p[\widetilde{N}]}~ := ~ \int_{{\cal F}}
\Spur\left(e^{-t\Delta_p[\widetilde{N}]}(x,x)\right)\,dx,\]
where $\Spur$ is the ordinary trace of an endomorphism of
a finite-dimensional vector space. \qed
\end{definition}

\begin{definition} \label{definition of determinant-class}
In the situation of Definition
\ref{definition of normalized trace} we
denote by $\Delta_p^{\perp}[\widetilde{N}]$ the operator
from the orthogonal complement of the kernel of
$\Delta_p[\widetilde{N}]$ to itself which is obtained from
$\Delta_p[\widetilde{N}]$ by restriction.
We call $N$ of {\em determinant-class} if for all $p \ge 0$
\[ \int_1^{\infty} \tr_{\Gamma}e^{-t\Delta_p^{\perp}[\widetilde{N}]}\,
\frac{dt}{t} ~
< ~ \infty. \qed\]
\end{definition}

The condition of determinant-class will be needed to define
$L^2$-analytic torsion and was introduced in
\cite[page 754]{Burghelea-Friedlander-Kappeler McDonald (1996a)}.
If all the Novikov-Shubin invariants of $N$ are positive
then $N$ is of determinant-class. There is the conjecture
that the Novikov-Shubin-invariants of a compact manifold
are always positive \cite{Novikov-Shubin(1986a)} or \cite[Conjecture
7.2]{Lott-Lueck (1995)}.
This has been verified for compact $3$-manifolds
whose prime factors with infinite fundamental groups
admit a geometric JSJT-decomposition
\cite[Theorem 0.1]{Lott-Lueck (1995)} and
for hyperbolic manifolds with or without boundary of finite volume
in \cite[section VII]{Lott  (1992a)}. If the fundamental group
is residually finite or amenable and $N$ is compact,
a proof that the manifold is of determinant-class is given in
\cite[Theorem A in Appendix A]{Burghelea-Friedlander-Kappeler (1996a)}
and \cite[Theorem 0.2]{Dodziuk-Mathai(1996a)}
using \cite[section 3]{Lueck (1994c)}.

\begin{lemma} \label{asymptotic expansion of the heat kernel}
In the situation of Definition
\ref{definition of normalized trace} the normalized trace
$\tr_{\Gamma}e^{-t\Delta_p[\widetilde{N}]}$
has for each $k \ge 0$  an asymptotic expansion for $t\to 0$
\[ \tr_{\Gamma}e^{-t\Delta_p[\widetilde{N}]} ~ = ~
\sum_{i=0}^k t^{-(m-i)/2}(a_i+b_i)+O(t^{(k-m+1)/2}). \]
Moreover, we have
\begin{eqnarray*}
a_i & = & \int_{{\cal F}} \alpha_i(x) \,dx;
\\
b_i & = & \int_{{\cal F} \cap \widetilde{\partial N}} \beta_i(x) \, dx,
\end{eqnarray*}
where ${\cal F}$ is a fundamental domain
for the $\Gamma$-action on $\widetilde{N}$,
such that ${\cal F} \cap \widetilde{N}$ is a fundamental domain
for the $\Gamma$-action on the preimage $\widetilde{\partial N}$
of $\partial N$ under the universal covering
$\widetilde{N} \longrightarrow N$,
$\alpha_i(x)$ is a density on $\widetilde{N}$ given
locally in terms of the metric
 and $\beta_i(x)$ is a density
on $\widetilde{\boundary N}$, which can be
computed locally out of the germ of
the metric at the boundary $\widetilde{\partial N}$. \qed
\end{lemma}
\proof
We begin with extending the result of Lott
\cite[Lemma 4]{Lott  (1992a)}) to manifolds with boundary. Namely
we want to prove the existence of constants $C_1,C_2 > 0$ independent of
$x \in \widetilde{N}$ and $t \in (0,\infty)$ such that for the covering 
projection
$q: \widetilde{N} \longrightarrow N$  we get
\begin{eqnarray}
\abs{e^{-t\Delta_p[\widetilde{N}]}(x,x) - e^{-t\Delta_p[N]}(q(x),q(x))}
\le C_1\cdot e^{-1/C_2t} \hspace{10mm}\forall x \in \widetilde{N}, t
\in(0,\infty). \label{eqn 1.111}
\end{eqnarray}
Fix a number $K > 0$ such that the restriction of
$q : \widetilde{N} \longrightarrow N$ to any ball $B_{2K} \subset N$
of radius $2K$  is a trivial covering. Consider $x \in \widetilde{N}$.
Choose a connected neighbourhood $V \subset N$ of $q(x)$
such that $B_K(x) \subset V \subset B_{2K}(x)$ and $V$
carries the structure of a Riemannian manifold for which the inclusion
of $V$ into $N$ is a smooth map respecting the Riemannian metrics.
We can find
$\widetilde{V} \subset \widetilde{N}$
such that $x\in\widetilde{V}$ and $\widetilde{V}$ carries the structure of a
Riemannian manifold for which the inclusion
of $\widetilde{V}$ into $\widetilde{N}$
is a smooth map respecting the Riemannian metrics
and $q$ restricted to $\widetilde{V}$ induces an isometric
diffeomorphism from $\widetilde{V}$ onto $V$. Since
Theorem \ref{kernComp} applies to $\widetilde{V} \subset \widetilde{N}$
and $V \subset N$, we obtain constants $C_1, C_2 > 0$ independent
of $x \in \widetilde{N}$ and $t \in (0,\infty)$
satisfying
\begin{eqnarray}
\abs{e^{-t\Delta_p[\widetilde{N}]}(x,x) -
e^{-t\Delta_p[\widetilde{V}]}(x,x)}
& \le & C_1/2\cdot e^{-1/C_2t}
\hspace{8mm}\forall t \in(0,\infty);
\label{eqn 1.112}
\\
e^{-t\Delta_p[\widetilde{V}]}(x,x) -
e^{-t\Delta_p[V]}(q(x),q(x)) & = & 0
\hspace{28mm} \forall t \in(0,\infty);
\label{eqn 1.113}
\\
\abs{e^{-t\Delta_p[N]}(q(x),q(x)) -
e^{-t\Delta_p[V]}(q(x),q(x))}
& \le &C_1/2\cdot e^{-1/C_2t}
\hspace{8mm}\forall t \in(0,\infty).
\label{eqn 1.114}
\end{eqnarray}
Now \eqref{eqn 1.111} follows from
\eqref{eqn 1.112}, \eqref{eqn 1.113} and \eqref{eqn 1.114}.

For each $k \ge 0$  there is an
asymtotic expansion for $t\to 0$
\[ \tr e^{-t\Delta_p[N]} ~ := \int_N \Spur_{\R}
\left(e^{-t\Delta_p[N]}(x,x)\right)\, dx
~ = ~
\sum_{i=0}^k t^{-(m-i)/2}(a_i[N]+b_i[N])+O(t^{(k-m+1)/2})\]
with
\begin{eqnarray*}
a_i[N] & = & \int_{N} \alpha_i[N](x) \,dx;
\\
b_i[N] & = & \int_{\partial N} \beta_i[N](x) \, dx,
\end{eqnarray*}
where $\alpha_i(x)[N]$ is a locally in terms of the metric given
density on $N$ and $\beta_i(x)$ is density
on $\boundary N$, which can be
computed locally out of the germ of
the metric at the boundary $\partial N$. This is a classical result
of Seeley \cite{Seeley (1969)} and Greiner  \cite{Greiner  (1971)}.
If we define the desired densities by
$\alpha_i  :=  \alpha_i[N] \circ q$ and
$\beta_i   :=  \beta_i[N] \circ q|_{\partial \widetilde{N}}$
the claim follows because $\vol(N) \cdot C_1\cdot e^{-1/C_2t} $ is
$O(t^{(N-m+1)/2})$. \qed

Now we recall the definition of $L^2$-analytic torsion
(here $\Gamma(s)=\int_{0}^{\infty} t^{s-1} e^{-t} ~ dt$).

\begin{definition}
\label{original definition of L^2-analytic torsion}
Let \mbox{$\widetilde{N} \longrightarrow N$} be the universal covering
of a compact connected Riemannian manifold $N$ of dimension $m$.
Suppose that $N$ is of determinant-class.
Define the {\em $L^2$-analytic torsion} of $N$ by
\begin{eqnarray*}
\anTor(N) & := & \sum_{p \ge 0} (-1)^p \cdot p \cdot\\
&  &  \hspace{10mm}  \left(\frac{d}{ds}\frac{1}{\Gamma(s)}
\int_0^{1} t^{s-1}\cdot
\tr_{\Gamma}e^{-t\Delta_p^{\perp}[\widetilde{N}]}\, dt\right|_{s=0}
+\left.\int_{1}^{\infty}
\tr_{\Gamma}e^{-t\Delta_p^{\perp}[\widetilde{N}]}\frac{dt}{t}\right).
\end{eqnarray*}
\end{definition}
Here the first integral is a priori only defined for $\Re(s) > m/2$
but Lemma \ref{asymptotic expansion of the heat kernel} ensures
that it has a meromorphic extension
to the complex plane with no pole in $s = 0$. The second integral
converges because of the assumption that $N$ is of determinant-class.

\begin{remark} \label{definitions make also sense in the
hyperbolic finite volume case}
Notice that Definition \ref{definition of normalized trace},
Definition \ref{definition of determinant-class},
Definition \ref{original definition of L^2-analytic torsion}
and Lemma \ref{new expression for zetafunction} carry over
to the case where $N$ is a not necessarily compact hyperbolic
complete Riemannian manifold with finite volume.
This follows from the fact that $\HH^m$ is homogeneous.\qed
\end{remark}

The proof of Theorem
\ref{Tor_converg} now splits into two separate parts:
we study the large time summand using
algebraically minded functional
analysis, namely, we refine the methods
of Lott-L\"uck \cite{Lott-Lueck  (1995)}.
The small time behaviour is handled
analytically via careful study of the heat
kernels. Here the main ingredient is the
``principle of not seeing the boundary''
due to Kac for functions, and
Dodziuk-Mathai \cite{Dodziuk-Mathai  (1996b)}
for forms. We give a short proof,
using only unit propagation speed, which was
suggested to us by Ulrich Bunke.

Finally we explain how Theorem \ref{TopTor=AnTor}
follows from Theorem \ref{Tor_converg}.
We use Notation \ref{notation for the ends}.
We want to apply the following result  of
Burghelea, Friedlander and Kappeler
\cite[Theorem 3.1]{Burghelea-Friedlander-Kappeler (1996a)}.

\begin{theorem} \label{comparision theorem for manifolds with boundary}
 Let $N$ be a compact Riemannian manifold of determinant-class
such that the Riemannian metric is a product near the boundary. Then
\[\anTor(N) ~ = ~ \topTor(N) + \frac{\ln(2)}{2} \cdot \chi(\boundary N).
\qed\]
\end{theorem}

 Therefore we equip $M_R$ with a new metric
which is a product near the boundary.
\begin{definition} \label{metric on M_R'}
Let $M'_R$ ($R\ge 1$) denote the same manifold as
$M_R$ but now equipped with a new Riemannian metric which is
equal to the old one outside $T_{R-1}\subset M_R$ and on $T_{R-1}$
given by
\[ g_R=du^2+e^{-2(\phi(R-u)u+(1-\phi(R-u))R)}dx^2, \]
where $\phi:[0,\infty] \to [0,1]$ is a smooth
function which is identically zero on $[0,\epsilon]$
and identically one on $[1-\epsilon,\infty)$ for some
$\epsilon > 0$.
 Similarly, let $E'_R$ be $E_R$ with the new Riemannian metric
 $du^2+e^{-2(\phi(u-R)R + (1-\phi(u-R))u)}dx^2$. \qed
\end{definition}

\begin{lemma}\label{isometries}
For arbitrary $r,R>0$, there is an isometric
diffeomorphism
\[ \tau_R: \widetilde{E_0} \longrightarrow \widetilde{E_R},\]
which induces by restriction isometric diffeomorphisms
$\widetilde{E_r} \to \widetilde{E_{r+R}}$,
$\widetilde{T_r} \to \widetilde{T_{r+R}}$,
$\widetilde{E_r'} \to \widetilde{E_{r+R}'}$ and
$\widetilde{T_r'} \to \widetilde{T_{r+R}'}$.
\end{lemma}
\proof
Define $\tau_R: [0,\infty) \times \R^{m-1} \longrightarrow
[R,\infty) \times \R^{m-1}$ by
$\tau_R(u,x) = (u+R,  e^R x)$. \qed

\begin{lemma} \label{variation of metrics at the ends}
For $m$ odd
\begin{eqnarray}
\lim_{R \to \infty} ~ \anTor(M_R) - \anTor(M'_R)  & = & 0; \label{vareq1}
\\
\topTor(M_R') - \topTor(\overline{M}) & = & 0
\hspace{6mm} \mbox{ for } R \ge 1.\label{vareq2}
\end{eqnarray}
\end{lemma}
\proof
Observe that  for $m$ odd $M$ is
$L^2$-acyclic by \cite{Dodziuk  (1979)}). The same holds
for $M_R$ by Corollary \ref{L2_acyclic} or
\cite[Theorem 1.1]{Cheeger-Gromov  (1985)}.
Let $\{g_u \mid u \in[0,1]\}$ be the obvious family of Riemannian metrics on
$M_R$ joining the hyperbolic metric on $M$ restricted
to $M_R$ with the metric $g_R$ on $M_R$ introduced in Definition
\ref{metric on M_R'}. Then we get from  Corollary \ref{metric anomaly}
 \[ \frac{d}{du}|_{u=0} \anTor(M_R,g_u) = \sum_p(-1)^{p}d_p[R,u],  \]
where $d_p[R,u]$ is an integral of a density $D_p[R,u]$
on $\boundary M_R$ which is given locally in
terms of the germ of the family of
metrics $g_u$ on $\boundary M_R$. If we pull back the density
$D_p[R,u]$ to a density $\widetilde{D_p[R,u]}$ on $\widetilde{\partial M_R}$
we can rewrite $d_p[R,u]$ as an integral over
the density $\widetilde{D_p[R,u]}$
over ${\cal F} \cap \widetilde{\partial M_R}$ for a fundamental domain
${\cal F}$. Notice that $\widetilde{D_p[R,u]}$ can be computed
locally in  terms of the germ of the family of
lifted metrics $\widetilde{g_u}$ on $\widetilde{\boundary M_R}$.
We can write $\widetilde{D_p[R,u]}$ as
$f_p[R,u] \cdot dvol_{\widetilde{\partial M_R}}$
for a function $f_p[R,u]$ on $\widetilde{\partial {M_R}}$.
Because the isometries in Lemma \ref{isometries} respect
the families of metrics, the functions
$f_p[R,u]$ are uniformly bounded in $R$ and $u$.
However, the volume of the domain of integration
 ${\cal F} \cap \widetilde{\partial M_R}$ tends to zero
if $R$ tends to $\infty$.
Hence \ref{vareq1} follows.

For \eqref{vareq2} we use the
fact that the topological $L^2$-torsion
is a homotopy invariant for
residually finite fundamental groups \cite[Theorem 0.5]{Lueck (1994c)}.
The isometries of the hyperbolic space $\HH^m$ can be
considered as a subgroup of $Gl(m+1,\R)$
\cite[Theorem A.2.4]{Benedetti-Petronio (1992)}.
Every finitely generated group
which possesses a faithful representation in
some $Gl(n,K)$ for any field $K$ is residually finite
\cite[Theorem 4.2]{Wehrfritz (1973)}.
We have seen that every  manifold which
possesses a complete hyperbolic metric with finite volume is homotopy
equivalent to a finite $CW$-complex.
It follows that  the
fundamental group of such a manifolds is
residually finite. \qed

Finally we can explain how Theorem \ref{TopTor=AnTor}
follows from Theorem \ref{Tor_converg}.
We begin with the case where $m = \dim(M)$ is odd.
Since $\boundary M$ is a union of tori
$\chi(\boundary M)$ is trivial.
Now one just applies Theorem
\ref{comparision theorem for manifolds with boundary}
for $N = M_R'$ and uses Lemma \ref{variation of metrics at the ends}.
Suppose that $m$ is even. Then the usual proof of Poincar\'e duality
for analytic torsion \cite[Theorem 2.3]{Ray-Singer(1971)}
extends directly to $L^2$-analytic torsion of $M$ (see
\cite[Proposition 16]{Lott  (1992a)}) and
shows $\anTor(M) = 0$. From \cite[Theorem 1.6,
Theorem 1.7 and Theorem 1.11]{Lueck (1994a)}
we conclude $\topTor(\overline{M}) = 0$. This finishes the proof
that Theorem \ref{Tor_converg} implies Theorem \ref{TopTor=AnTor}.

The rest of this paper is devoted to the proof of Theorem
\ref{Tor_converg} (and also of Theorem \ref{kernComp} and
Corollary \ref{metric anomaly} which we have already used above).



\typeout{--------------------- section 2-----------------------}

\section{Analysis of the heat kernel}
\label{Analysis of the heat kernel}
\subsection{Standard Sobolev estimates}

\begin{definition}
Let  $N$ be a Riemannian manifold with boundary $x\in\boundary N$.
We define
the {\em boundary exponential map}
\[ \exp_x^\boundary:[0,\infty)\times T_x\boundary N \to N: (u,y)
\mapsto  \nu(u,\exp_x^{\boundary N}(y)).\]
Here $\exp^{\boundary N}$ is the exponential map of the boundary with its
induced Riemannian metric and $\nu$ denotes the geodesic flow of the inward
unit normal field.\par

If we identify $T_x^{\boundary N}$ with $\R^{m-1}$ via an orthonormal frame
and restrict the boundary exponential map to a subset where it is a
diffeomorphism onto its image, we get so called {\em boundary normal
coordinates}.

For an interior point the coordinates
induced from the exponential map are
called {\em Gaussian coordinates}.
The term {\em normal coordinates} is used
to denote Gaussian coordinates as well as boundary normal coordinates.
\qed
\end{definition}

\begin{definition}\label{bc}
  A form $\omega$ on a Riemannian manifold with boundary fulfills
{\em  relative boundary conditions for $\Delta^m$, $m\ge 1$}, if
\[ i^*(\Delta^j\omega)=0\quad\text{and}\quad i^*(\Delta^j\delta\omega)=0
\qquad\text{for $0\le j\le m-1$}. \]
It fulfills {\em absolute boundary conditions for $\Delta^m$} if
\[ i^*(*\Delta^j\omega)=0\quad\text{and}\quad i^*(*\Delta^jd\omega)=0
\qquad\text{for $0\le j\le m-1$}. \]
Here and in the sequel, $i:\boundary M\to M$ always denotes the inclusion of
the boundary. \qed

Remember that both boundary value problems are elliptic.
The next result is a standard elliptic estimate
(compare f.i.~\cite[1.24]{Cheeger-Gromov-Taylor (1982)}
 if $\boundary M=\emptyset$). Since we do not know an explicite
reference for manifolds with boundary
we include a proof here (compare also \cite[Ch.5, Sec.9]{Taylor (1996a)}).
\end{definition}

\begin{notation}
For real numbers $a,b,c>0$ we abbreviate
\[ a\stackrel{c}{\le}b\qquad\text{for}\quad a\le cb \]
\end{notation}

\begin{theorem}\label{sob_est}
Let $M$ be a Riemannian manifold with boundary
of dimension $m$, $x_0\in M$ and
$r$ sufficiently small so that
$B_r(x_0)\subset M$ is diffeomorphic in normal coordinates to
$B_r(x)\subset \R^m_{\ge 0}$ for some
$x\in \R^m_{\ge 0}$, where $\R^m_{\ge 0}$ denotes the half space.
(The location of $x$ depends on the question whether $B_r(x_0)$ meets the
boundary of $M$.)\par

Then we can find $C>0$ so that for all $\omega\in C^\infty$
 which fulfill either absolute or
relative boundary conditions
\[ \abs{\omega(x_0)}^2
\le
C\sum_{i=0}^m \abs{\Delta^i \omega}^2_{L^2(B_r(x_0))}. \]
The constant $C$ is a smooth function of the coefficients of the
Riemannian metric, its inverse and their
derivatives in normal coordinates of
$B_r(x_0)$. Moreover, it depends on $r$ (it becomes larger if
$r$ becomes smaller).
\end{theorem}
\begin{proof}
We start with the following formula:

\begin{lemma}\label{formel}
For every $k$-form $f$ which fulfills either absolute or relative boundary
conditions for $\Delta$ and for every function $\phi$
\[ \abs{d(\phi f)}^2_{L^2(M)}+\abs{\delta(\phi f)}^2_{L^2(M)}
~ = ~
(\Delta f,\phi^2 f)_{L^2(M)} + \abs{f\wedge d\phi}^2_{L^2(M)}
 + \abs{*f\wedge d\phi}^2_{L^2(M)}. \]
\end{lemma}
\begin{proof}

\begin{eqnarray}
(d(\phi f),d(\phi f))
& = &
(d\phi\wedge f,d(\phi f))+(\phi df,d(\phi f)) \nonumber
\\
& = &
(d\phi\wedge f,d\phi\wedge f)+(d\phi\wedge f,\phi df) +
 (df,\phi d(\phi f))\nonumber
\\
& = &
(d\phi\wedge f,d\phi\wedge f)+(d\phi\wedge f,\phi df) +
(df,d(\phi^2 f)- d\phi\wedge(\phi f)) \nonumber
\\
& = &
\abs{d\phi\wedge f}^2 +(\delta d f,\phi^2 f).
\label{eqn 2.800}
\end{eqnarray}
For the last equation we have used the fact
that the boundary contribution of the integration by parts
vanishes if either $i^* f=0$ or $i^*(*df)=0$ and that
the first equation holds if $f$ fulfills relative
boundary conditions, and the second equation holds if $f$
satisfies absolute boundary conditions.
\begin{eqnarray}
(\delta(\phi f),\delta(\phi f))
& = &
(d*(\phi f),d*(\phi f)) \nonumber
\\
& = &
(d(\phi *f),  d(\phi *f)) \nonumber
\\
& = &
\abs{d\phi\wedge(*f)}^2 +(\delta d(*f),\phi^2 *f) \nonumber
\\
& = &
\abs{d\phi\wedge(*f)}^2 +
 ( (-1)^{mk+m+1} d*d* f,\phi^2 f) \nonumber
\\
& = &
\abs{d\phi\wedge(*f)}^2 +(d\delta f,\phi^2 f).
\label{eqn 2.801}
\end{eqnarray}
We get the third equation above by applying \eqref{eqn 2.800}.
Now add \eqref{eqn 2.800} and \eqref{eqn 2.801}
\end{proof}

The next lemma will be needed several times.
In the sequel $\abs{\cdot }^2_{H^k(\R^m_{\ge 0})}$
denotes the Sobolev norm on $\R^m_{\ge 0}$ with the standard
Euclidean metric, and a function $f$ with compact support in an open
subset $V \subset \R^m_{\ge 0}$ is extended to $\R^m_{\ge 0}$
by zero on the complement of $V$. Moreover, $\Delta$ stands
for the Laplacian on the Riemannian manifold $M$.

\begin{lemma}\label{kl_lemma}
In the situation of Theorem \ref{sob_est}
let $\phi,\psi$ be  functions with
$\supp\phi\subset\supp\psi\subset B_r(x_0)$ and $\psi\equiv 1$ on $\supp\phi$.
We use normal coordinates to identify $B(r,x_0)$ with a
subset of $\R^m_{\ge 0}$.
Then for $k,l\ge 0$ and every form $\omega$
which fulfills either absolute or relative boundary conditions for
$\Delta^{l}$
\begin{equation}
  \label{kl_est}
  \abs{\phi\omega}^2_{H^{k+2l}(\R^m_{\ge 0})} \stackrel{C}{\le}
 \abs{\phi\omega}^2_{H^k(\R^m_{\ge 0})}
+ \abs{\phi\Delta^l\omega}^2_{H^k(\R^m_{\ge 0})} +
 \abs{\psi \omega}^2_{H^{k+2l-1}(\R^m_{\ge 0})}.
\end{equation}
The constant $C$ smoothly depends  on the coefficients of the
 metric, its inverse
and their derivatives in normal coordinates in $B(r,x_0)$.
 In addition, it depends on $\phi$.
\end{lemma}
\begin{proof}
We prove only the case of relative boundary conditions.
\begin{eqnarray}
\abs{\phi\omega}^2_{H^{2l+k}(\R^m_{\ge 0})}
 & \stackrel{C_1}{\le} &
 \abs{\Delta^{l}(\phi\omega)}_{H^k(\R^m_{\ge 0})}^2
+ \abs{\phi\omega}^2_{H^k(\R^m_{\ge 0})}
+ \sum_{j=0}^{l-1}
\abs{i^*(\Delta^j(\phi \omega))}^2_{H^{2l+k-2j-1/2}(\R^{m-1})}
\nonumber
\\
& & \hspace{40mm}
+ ~ \abs{i^*(\Delta^j\delta(\phi \omega))}^2_{H^{2l+k-2j-3/2}(\R^{m-1})}
\label{eqn 2.802}
 \\
 & \stackrel{2}{\le} &
\abs{\phi\omega}^2_{H^k(\R^m_{\ge 0})} +
 \abs{\phi\Delta^l\omega}_{H^k(\R^m_{\ge 0})}^2  +
 \abs{[\Delta^l,\phi]\omega}^2_{H^k(\R^m_{\ge 0})}
\nonumber
 \\
& &
+ ~ \sum_{j=0}^{l-1}
\abs{\phi\underbrace{i^*(\Delta^j\omega)}_{=0}}^2_{H^{2l+k-2j-1/2}(\R^{m-1})}
+ \abs{i^*([\Delta^j,\phi] \omega)}^2_{H^{2l+k-2j-1/2}(\R^{m-1})}
\nonumber
 \\
& &
+ ~ \abs{\phi\underbrace{i^*(\Delta^j\delta
\omega)}_{=0}}^2_{H^{2l+k-2j-3/2}(\R^{m-1})}
+ \abs{i^*([\Delta^j\delta,\phi] \omega)
 }^2_{H^{2l+k-2j-3/2}(\R^{m-1})} \label{eqn 2.803}
 \\
& \stackrel{C_2}{\le} &
\abs{\phi\omega}^2_{H^k(\R^m_{\ge 0})} +
 \abs{\phi\Delta^l\omega}_{H^k(\R^m_{\ge 0})} +
 \abs{\psi\omega}^2_{H^{2l+k-1}(\R^m_{\ge 0})}.
 \label{eqn 2.804}
\end{eqnarray}
For \eqref{eqn 2.802} we use the fact that $\Delta$ with relative boundary
conditions
is elliptic. Then, the same is true for its $l^{th}$ power and
therefore we can estimate the Sobolev norm as indicated
(compare for instance \cite[4.15]{Schick (1996)}).
Note that the constant $C_1$ (smoothly) depends 
on the coefficients of the
boundary value problem, which are determined by
the Riemannian metric in $B(r,x_0)$
and its derivatives (again in a smooth way). Therefore,
$C_1$ smoothly depends  on the 
 metric as desired.
\par

For \eqref{eqn 2.803} we use the fact that
$\omega$ fulfills relative boundary conditions
and we denote by $[\Delta^j,\phi]$ the commutator of $\Delta^j$
and the operator
given by multiplication with $\phi$. \par

In \eqref{eqn 2.804} we use the fact that the commutator $[\Delta^l,\phi]$
 is a differential operator of order $2l-1$.
Similarly, the commutators $[\Delta^j,\phi]$ have order $2j-1$ and
$[\Delta^j\delta,\phi]$ have order $2j$. Also, $i^*$ is the
composition of the trace map
$H^s(\R^m_{\ge 0})\to H^{s-1/2}(\R^{m-1})$ for $s>1/2$ which simply
restricts to the boundary but leaves the bundle unchanged, and the 
operator of order $0$ (bounded on $H^s$ for all $s$) which projects
onto the ``tangential'' components of forms. Therefore
the statement follows.
Notice that we can replace $\omega$ by $\psi\omega$ above in the last
summand as $i^*([\Delta^j,\phi] \omega) = i^*([\Delta^j,\phi] \psi\omega)$.
It remains to check what the constant $C_2$ depends on,
i.e. which norm the operators mentioned above have.
They all are differential operators and their norm depends on the
coefficients. These coefficients are determined by the coefficients of
$\Delta$ and $\delta$ and $i^*$
and their derivatives, but now also by $\phi$ and its
derivatives.
So, all together, we have the (smooth) dependence on $\phi$
 and  the Riemannian metric as desired.
\end{proof}

Inductive application of the following lemma will prove
Theorem \ref{sob_est}.
\begin{lemma}\label{2step}
Adopt the situation of Lemma \ref{kl_lemma} and assume that $\phi$ is a
cutoff function, i.e. $\phi:M\to[0,1]$. Then for every $k\ge 0$ and
every $\omega$
which fulfills either absolute or relative boundary conditions we get
  \begin{equation}\begin{split}
    \label{2step_est}
    \abs{\phi\omega}^2_{H^{2k}(\R^m_{\ge 0})} \stackrel{C}{\le}&
 \abs{\Delta^k\omega}^2_{L^2(B(r,x_0))} +
 \abs{\Delta^{k-1}\omega}^2_{L^2(B(r,x_0))} +
 \abs{\Delta\omega}^2_{L^2(B(r,x_0))} +
 \abs{\omega}^2_{L^2(B(r,x_0))}\\
 & +
 \abs{\psi\omega}^2_{H^{2k-2}(\R^m_{\ge 0})}.
\end{split}\end{equation}
On $L^2(B(r,x_0))$ we use the norm induced from the Riemannian metric on $M$.
$C$ smoothly depends  on the coefficients of the Riemannian metric tensor and
its inverse in normal coordinates and on their derivatives.
Also it smoothly depends  on $\phi$ and $\psi$.
If $k=1$ we can replace the norm $\abs{\cdot}_{H^{2k-2}(\R^m_{\ge 0})}$
by $\abs{\cdot}_{L^2(B(r,x_0))}$ and the statement remains true.
\end{lemma}
\begin{proof}
We proof only the case of relative boundary conditions.\par

Since we identified $B(r,x_0)$ with a subset of $\R^m$ we have two norms
on $L^2(B(r,x_0))$: the one used in the statement of
Lemma \ref{2step} coming from the
Riemannian metric on $M$ and the one coming from Euclidian $\R^m$.
 We find a constant $C$ depending smoothly on the
Riemannian metric tensor and its inverse so that
\begin{eqnarray}
\abs{\cdot}_{L^2(\R^m)}
& \stackrel{1/C}{\le} &
\abs{\cdot}_{L^2(B(r,x_0))}
~ \stackrel{C}{\le} ~ \abs{\cdot}_{L^2(\R^m)}. \label{eqn 2.90}
\end{eqnarray}
In particular the last statement in Theorem \ref{2step}
for $k=1$ follows from the rest
and \eqref{eqn 2.90}.\par

Choose a cutoff function $\phi_1:M\to [0,1]$
 with $\phi_1\equiv 1$ on $\supp\phi$ and
$\psi\equiv 1$ on $\supp\phi_1$.
\begin{eqnarray}
\abs{\phi\omega}^2_{H^{2k}(\R^m_{\ge 0})}
&\stackrel{C_1}{\le} &
 \abs{\phi\omega}^2_{L^2(B(r,x_0))} +
\abs{\phi\Delta^k\omega}^2_{L^2(B(r,x_0))}
 +\abs{\phi_1\omega}^2_{H^{2k-1}(\R^m_{\ge 0})} \label{eqn 2.91}
\\
& \stackrel{C_2}{\le} &
\abs{\omega}_{L^2(B(r,x_0))}^2 + \abs{\Delta^k\omega}_{L^2(B(r,x_0))}^2
+ \abs{\phi_1\omega}^2_{H^1(\R^m_{\ge 0})}  \nonumber
 \\ & & \hspace{30mm}
 + \abs{\phi_1\Delta^{k-1}\omega}^2_{H^1(\R^m_{\ge 0})} +
 \abs{\psi\omega}^2_{H^{k-2}(\R^m_{\ge 0})}. \label{eqn 2.92}
\end{eqnarray}
Inequality \eqref{eqn 2.91} follows from
\eqref{kl_est} with $k=0$ and $l$ replaced by $k$, and
\eqref{eqn 2.90}. The constant $C_1$ depends only on $\phi$ and the
Riemannian metric in $B(r,x_0)$.
We conclude \eqref{eqn 2.92}
from $0\le\phi\le 1$ and \eqref{kl_est} now with $k=1$ and
$l=k-1$. We clearly can choose $\phi_1$ depending smoothly on $\phi$ and
$\psi$ so that $C_2$ depends smoothly on $\phi$,
the coefficients of the Riemannian metric tensor
and its inverse in normal coordinates
and on their derivatives.\par

It remains to estimate the two $H^1$-summands appearing in \eqref{eqn 2.92}.
\begin{eqnarray}
\abs{\phi_1\omega}^2_{H^1(\R^m_{\ge 0})}
& \stackrel{C_3}{\le} &
 \abs{\phi_1\omega}^2_{L^2(B(r,x_0))} +
 \abs{d(\phi_1\omega)}^2_{L^2(B(r,x_0))}
 +\abs{\delta(\phi_1\omega)}^2_{L^2(B(r,x_0))} \nonumber
 \\
& & \hspace{40mm}
 +~ \abs{\phi_1 \underbrace{i^*\omega}_{=0}}^2_{H^{1/2}(\R^{m-1})}
 \label{eqn 2.70}
\\
& \le &
\abs{\omega}^2_{L^2(B(r,x_0))}
+ (\Delta\omega,\phi_1^2\omega)_{L^2(M)} +
 \abs{\omega\wedge d\phi_1}^2_{L^2(M)} \nonumber
 \\
& & \hspace{40mm}
+ ~ \abs{*\omega\wedge d\phi_1}^2_{L^2(M)} \label{eqn 2.71}
\\
& \stackrel{C_4}{\le} &
\abs{\omega}^2_{L^2(B(r,x_0))} +
\abs{\Delta\omega}_{L^2(B(r,x_0))}\cdot\abs{\phi_1^2\omega}_{L^2(B(r,x_0))}
\label{eqn 2.72}
\\
& \le &
\abs{\omega}^2_{L^2(B(r,x_0))} +
\abs{\Delta\omega}_{L^2(B(r,x_0))}^2 +
\abs{\phi_1^2\omega}_{L^2(B(r,x_0))}^2 \label{eqn 2.73}
\\
& \le &
2 \cdot \abs{\omega}^2_{L^2(B(r,x_0))} +
\abs{\Delta\omega}_{L^2(B(r,x_0))}^2.  \label{eqn 2.74}
\end{eqnarray}
We get \eqref{eqn 2.70}
since $(d,\delta; i^*)$ is an elliptic boundary value problem
and $\omega$ fulfills relative boundary conditions. Obviously, $C_3$ is
determined in the same way as above. We conclude
\eqref{eqn 2.71} from Lemma \ref{formel} and $0\le\phi_1\le 1$.
Equation \eqref{eqn 2.72} follows from  the Cauchy-Schwarz inequality.
The constant $C_4$ in \eqref{eqn 2.72} involves
$\sup_{x\in B(r,x_0)}\abs{d\phi_1(x)}$).
Equation \eqref{eqn 2.73}
follows from $ab\le a^2+b^2$ for arbitrary real numbers $a$ and $b$.
We get \eqref{eqn 2.74} from $0\le\phi_1\le 1$.
An identical argument shows for the second $H^1$-summand
\begin{eqnarray}
\abs{\phi_1\Delta^{k-1}\omega}^2_{H^1(\R^m_{\ge 0})}
& \stackrel{c_5}{\le}&
2 \cdot \abs{\Delta^{k-1}\omega}^2_{L^2(B(r,x_0))} +
\abs{\Delta^k\omega}_{L^2(B(r,x_0))}^2.  \label{eqn 2.75}
\end{eqnarray}
Now Lemma \ref{2step} follows from
\eqref{eqn 2.91} to \eqref{eqn 2.75}.
\end{proof}

Finally we give the proof of
Theorem \ref{sob_est}.
Choose a cutoff function $\phi_1:M\to [0,1]$ with $\supp\phi_1\subset
B(r,x_0)$ and so that $\phi_1\equiv 1$ in a neighbourhood of $x_0$. By
Sobolev's lemma
\[\abs{\omega(x_0)}^2 = \abs{\phi_1\omega(x_0)}^2
  \stackrel{C_{2m}}{\le} \abs{\phi_1\omega}^2_{H^{2m}(\R^m_{\ge 0})}. \]
This is a computation entirely in $\R^m_{\ge 0}$, therefore
 $C_{2m}$ depends only
on the dimension. Theorem \ref{sob_est}
follows now by  an inductive application of Lemma \ref{2step}.
To do this, we have to choose a sequence of cutoff functions $\phi_i:M\to
[0,1]$ with $\supp \phi_i\subset B(r,x_0)$ for all $i$ and so that
$\phi_{i+1}\equiv 1$ on $\supp \phi_i$.
Clearly, we can construct this sequence
depending only on $r$. Since the
derivatives of these functions become larger
if $r$ becomes smaller, the constant
$C$ of Theorem \ref{sob_est} has to become larger,
too.
\end{proof}

\subsection{Comparison of heat kernels}
\label{heat_ker_comp}

In this section we use unit propagation speed (as in
\cite{Cheeger-Gromov-Taylor (1982)}) to prove the ``principle of
not feeling the boundary'' of M.\ Kac. Similar results have been proved by
Dodziuk-Mathai \cite{Dodziuk-Mathai (1996b)} with a  more complicated
argument, involving finite propagation speed and Duhamel's principle.
Moreover, their method does not yield the statement in the
generality we prove (and need) it.
The method we use was suggested
to us by Ulrich Bunke during the meeting on
Dirac operators 1997 at the Banach Center in Warzawa and uses
ideas of \cite{Dodziuk-Mathai (1996b)}.

The next definition extends the notion of a Riemannian manifold of bounded
geometry to manifolds with boundary (compare Schick
\cite[chapter 3]{Schick (1996)}, where these manifolds
are discussed in greater detail).
\begin{definition}
  A Riemannian manifold $(N,g)$ (possibly with boundary) is called a manifold
of {\em bounded geometry} if constants $C_k;\;k\in\N$ and
 $R_I, R_C>0$ exist, so that the following holds:
 \begin{enumerate}
\item the geodesic flow of the unit inward normal field induces a
diffeomorphism of $[0,2R_C)\times \boundary N$ onto its image $C(\boundary N)$,
 the geodesic
collar. Let $\pi: C(\boundary N)\to \boundary N$ be the corresponding
projection;

\item
$\forall x \in N$ with $d(x,\boundary N)>R_C/2$
 the exponential map $T_xN\to N$ is
a diffeomorphism on $B_{R_I}(0)$;

\item $\forall x\in N$ with $d(x,\boundary N)<R_C$ we have boundary normal
coordinates defined on $[0,R_C)\times
(B(0,R_C)\subset T_{\pi(x)}^{\boundary N})$;

\item For every $k\in\N$ and every $x\in\N$ the derivatives up to order $k$ of
the Riemannian metric tensor $g_{ij}$ and its inverse $g^{ij}$ in Gaussian
coordinates (if $d(x,\boundary N)>R_C/2$) or in normal boundary coordinates
(if $d(x,\boundary N)<R_C$) resp., are bounded by $C_k$. \qed
 \end{enumerate}
\end{definition}

\begin{example}
  Any covering of a compact manifold with the induced metric is a manifold of
bounded geometry. \qed
\end{example}

\begin{theorem}\label{kernComp}
Let $N$ be a Riemannian manifold possibly with boundary
which is of bounded geometry.
Let $V\subset N$ be a closed subset
which carries the structure of a Riemannian manifold
of the same dimension as $N$ such that the inclusion of $V$ into $N$
is a smooth map respecting the Riemannian metrics.
(We make no assumptions about the
boundaries of $N$ and $V$ and how they intersect.)
For fixed $p \ge 0$ let $\Delta[V]$ and $\Delta[N]$ be
the Laplacians on $p$-forms on $V$ and $N$, considered as
unbounded operators with either
absolute boundary conditions
or with relative boundary conditions (see Definition \ref{bc}).
Let $\Delta[V]^ke^{-t\Delta[V]}(x,y)$ and $\Delta[N]^k e^{-t\Delta[N]}(x,y)$
be the corresponding smooth integral kernels. Let $k$ be a non-negative
integer.

Then there is a monotone decreasing function
$C_k(K)$ from $(0,\infty)$ to $(0,\infty)$
which depends only on the geometry of
$N$ (but not on $V$, $x$, $y$, $t$) and a constant $C_2$ depending
only on the dimension of $N$ such that for all
$K > 0$ and $x,y\in V$ with
$d_V(x):=d(x,N-V)\ge K$, $d_V(y)\ge K$ and all $t > 0$:
\[ \abs{\Delta[V]^k e^{-t\Delta[V]}(x,y)-\Delta[N]^k e^{-t\Delta[N]}(x,y)}\le
C_k(K) e^{-(d_V(x)^2+d_V(y)^2+d(x,y)^2)/C_2t}. \]
\end{theorem}
\proof In the sequel $\Delta$ stands for both
$\Delta[N]$ and $\Delta[V]$. The Fourier transform of an $L^1$-function
$f: \R \longrightarrow \C$ is defined by
$\widehat{f}(\xi) = 1/\sqrt{2\pi} \cdot
\int_{-\infty}^{\infty}f(s) e^{-is\xi} ds$. Because
of the well-known rules
\begin{eqnarray*}
\widehat{e^{-s^2/2}} & = & e^{-\xi^2/2};
\\
\widehat{\frac{d^mf}{ds^m}} & = & i^m \xi^m \cdot \widehat{f};
\\
\lambda \cdot \widehat{g}(\lambda\xi) & = & \widehat{f}(\xi) \hspace{10mm}
\mbox{ for } g(s) := f(\lambda s),
\end{eqnarray*}
we conclude for $\xi \in \R$
\begin{eqnarray*}
\frac{(-1)^m}{\sqrt{\pi t}} \cdot
\int_0^{\infty} \left(\frac{d^{2m}}{ds^{2m}}e^{-s^2/4t}\right)
 \cos(s\xi)~ ds
 =&
\frac{(-1)^m}{\sqrt{2t}} \cdot
\Re\left(\left(\widehat{
\frac{d^{2m}}{ds^{2m}}e^{-s^2/4t}}\right)(\xi)\right)\\
= & \xi^{2m} e^{-t\xi^2}.
\end{eqnarray*}
By the spectral theorem applied to $\sqrt{\Delta}$ we get for
non-negative integers $l$, $m$ and $k$
\begin{eqnarray}
\Delta^m\Delta^l \Delta^k e^{-t\Delta}   & = &
\frac{(-1)^{l+m+k}}{\sqrt{\pi t}} \cdot
\int_0^{\infty} \frac{d^{2(m+l+k)}}{ds^{2(m+l+k)}}e^{-s^2/4t}
\cos(s\sqrt{\Delta})~ ds.
\label{eqn 2.1}
\end{eqnarray}
Choose $x_0,y_0\in V$ with $d_V(x_0)$, $d_V(y_0)\ge K$.
Notice for the sequel that the ball with radius $R$ less or equal to $K$
around $y_0$ in $N$ and the one in $V$ agree and will be denoted
by $B_R(y_0)$. Moreover, the intersection of $B_R(y_0)$ with $\partial N$
and with $\partial V$ agree. For a smooth
$p$-form $u$ with $\supp(u)\subset B_{K/4}(y_0)$
which satisfies the absolute or relative boundary conditions and hence
lies in the domain for both $\Delta[V]$ and $\Delta[N]$, we consider
now the function $f$ on $V$ given by
\[f ~ := ~ \left(\Delta[N]^k e^{-t\Delta[N]}-\Delta[V]^k
 e^{-t\Delta[V]}\right)u. \]
We conclude from \eqref{eqn 2.1}
\begin{eqnarray}
\lefteqn{\Delta^m \Delta^l f} \nonumber
\\
& = &
\int_0^\infty
 t^{-2(m+l+k)-1/2}P_{m,l,k}(s,\sqrt{t})
e^{-s^2/4t}(\cos(s\sqrt{\Delta[N]})-\cos(s\sqrt{\Delta[V]}))u\,ds~,
\label{f_Lap}
\end{eqnarray}
where $P_{m,l,k}$ is a universal polynomial with real coefficients.
Next we show
\begin{eqnarray}
\cos(s\sqrt{\Delta[N]})u  -  \cos(s\sqrt{\Delta[V]})u & = & 0
\hspace{5mm} \mbox{ on }B_{K/4}(x_0) \nonumber
\\ & &
\hspace{5mm} \mbox{for } s \le \max\{d_V(y_0)/2,d(x_0,y_0)/2\}.
\label{eqn 2.2}
\end{eqnarray}
Note that $\cos(s\sqrt{\Delta})u$ fullfils
the wave equation with initial data $u$.
By unit propagation speed for the wave
equation on manifolds with boundary
\cite[6.1]{Taylor (1996a)}),
$\supp(\cos(s\sqrt{\Delta})u)\subset B_{d_V(y_0)}(y_0)\subset V$ for
$s\le d_V(y_0)/2<d_V(y_0)-K/4$. Moreover, because of $\supp u \subset
B_{K/4}(y_0)$ and the uniqueness of solutions of the
wave equation we get on $V$
\begin{eqnarray}
\cos(s\sqrt{\Delta[N]})u  & = &
\cos(s\sqrt{\Delta[V]})u \quad\text{for }s\le
d_V(y_0)/2.
\label{eqn 2.3}
\end{eqnarray}
If $d(x_0,y_0)\ge d_V(y_0)\ge K$,
we know from unit propagation speed that
\begin{eqnarray}
\supp(\cos(s\sqrt{\Delta})u)\cap B_{K/4}(x_0) & = & \emptyset
\qquad\text{for }
 s\le d(x_0,y_0)/2 < d(x_0,y_0)-K/2,
\label{eqn 2.4}
\end{eqnarray}
since $\supp u \subset B_{K/4}(y_0)$.
Now \eqref{eqn 2.2} follows from \eqref{eqn 2.3} and \eqref{eqn 2.4}.
Since $\abs{\cos(s\xi)}\le 1$ and
hence $\abs{\cos(s\sqrt{\Delta})u}_{L^2}\le\abs{u}_{L^2}$
we conclude from \eqref{f_Lap} and \eqref{eqn 2.2}
\[\begin{split}
 \abs{\Delta^m \Delta^lf}_{L^2(B_{K/4}(x_0))} &\le
2\left(\int_{ \max\{d_V(y_0)/2,d(x_0,y_0)/2\}}^\infty
t^{-(2(m+l+k)+1/2)}P_{m,l,k}(s,\sqrt{t})e^{-s^2/4t}\, ds\right)
\abs{u}_{L^2}\\
&\le C_{m,l,k}e^{- \max\{d_V(y_0)/2,d(x_0,y_0)/2\}^2/C_{m,l}t} \abs{u}_{L^2}
\end{split}\]
by an elementary estimate of the integral. Since $N$ is of bounded geometry
and the heat kernel fulfills absolute boundary conditions, the
 elliptic estimates of Theorem \ref{sob_est} yield
pointwise bounds
\begin{eqnarray}
\abs{\Delta^lf(x_0)} & \le & D_l(K) \cdot
e^{- \max\{d_V(y_0)/2,d(x_0,y_0)/2\}^2/E_l t} \abs{u}_{L^2(B_{K/4}(x_o))},
\label{eqn 2.5}
\end{eqnarray}
where $D_l(K)$ is a monotone decreasing function in $K > 0$
which is given in universal expressions involving
the norm of curvature and a bounded number of its derivatives
on $B_{K/4}(x_0)$ and is independent of $x_0$, $y_0$ and $t$, and $E_l > 0$
depends only on the dimension of $N$.

Now $\Delta^lf(x_0)=\int
\Delta_y^l(\Delta^ke^{-t\Delta[N]}(x_0,y)-\Delta^k
e^{-t\Delta[V]}(x_0,y))u(y)~dy$.
Since the estimates \eqref{eqn 2.5}
hold for a dense subset $\{u\}\subset L^2(B_{K/4}(y_0))$ we conclude
\[
\abs{\Delta_y^l(\Delta^ke^{-t\Delta[N]}
(x_0,y)-\Delta^k e^{-t\Delta[V]}(x_0,y))}_{L^2(B_{K/4}(y_
0))} \le D_{l,k}(K)\cdot
e^{-\max\{d_V(y_0)/2,d(x_0,y_0)/2\}^2/E_l t}. \]
The very same reasoning as above yields pointwise bounds
\begin{eqnarray}
\abs{(\Delta^k e^{-t\Delta[V]}(x_0,y_0)-\Delta^k
e^{-t\Delta[N]}(x_0,y_0)} & \le &
C_k(K) e^{-\max\{d_V(y_0)/2,d(x_0,y_0)/2\}^2/C_2 t},
 \label{eqn 2.6}
\end{eqnarray}
where $C_k(K) > 0$ is a monotone decreasing function in $K > 0$
which is independent of $x_0$, $y_0$ and $t$,  and $C_2 > 0$
depends only on the dimension of $N$.
Since the heat kernel is symmetric we also have
\begin{eqnarray}
\abs{(\Delta^k e^{-t\Delta[V]}(x_0,y_0)-\Delta^k
e^{-t\Delta[N]}(x_0,y_0)} & \le & C_k(K)
e^{-\max\{d_V(x_0)/2,d(x_0,y_0)/2\}^2/C_2 t}.
 \label{eqn 2.7}
\end{eqnarray}
 Theorem \ref{kernComp}
follows from \eqref{eqn 2.6} and \eqref{eqn 2.7}.\qed

\begin{theorem}\label{ker_est}
Let $N$ be a Riemannian manifold possibly with boundary
which is of bounded geometry.
For $t_0>0$ we find $c(t_0)$,
depending on the bounds of the metric tensor and its inverse
and finitely many of their derivatives in normal coordinates
 (but not on $x, y\in N$) so
that
\[ \abs{ e^{-t\Delta[V]}(x,y)}\le c(t_0)\qquad\forall t\ge t_0. \]
\end{theorem}
\begin{proof}
  The norm of the bounded operator
  $\Delta^m \Delta^l e^{-t\Delta}$ is, by the
spectral theorem, bounded by
$\sup_{x\ge 0} x^{m+l}e^{-tx}$. For $t\ge t_0$  this
is bounded by some
 constant $C_{t_0}$. This $L^2$-bound can
be converted to a pointwise bound
exactly in the same way as above. Here, via
the Sobolev estimates of Theorem \ref{sob_est},
 the constant depends on the geometry of $N$.
\end{proof}

\subsection{Convergence of the small $t$ part of determinants}
\label{Convergence of the small $t$ part of determinants}
In this section, we study the small
$t$ summands in the Definition
\ref{original definition of L^2-analytic torsion} of
$\anTor(N)$. We use Lemma
\ref{asymptotic expansion of the heat kernel}
to rewrite the first integral involving small $t$ in
Definition \ref{original definition of L^2-analytic torsion}
in a form which does not involve meromorphic extension.
Namely, for
\begin{eqnarray*}
d^{sm}_p & := & \int_0^1\left( \tr_\Gamma
e^{-t\Delta_p^{\perp}[\widetilde{N}]}-
\sum_{i=0}^m t^{-(m-i)/2}(a_i+b_i)\right)
\,\frac{dt}{t} + \sum_{i=0}^m c(i,m)(a_i+b_i);
\\
c(i,m) & := & - \frac{m-i}{2} \hspace{5mm} \mbox{ for }  i \ne m;
\\
c(m,m) & := & - \left.\frac{d\Gamma}{ds}\right|_{s=1},
\end{eqnarray*}
we want to show
\begin{lemma} \label{new expression for zetafunction}
\[
\left.\frac{d}{ds}\frac{1}{\Gamma(s)}
\int_0^{1} t^{s-1}\cdot
\tr_{\Gamma}e^{-t\Delta_p^{\perp}[\widetilde{N}]}\, dt\right|_{s=0}
~ = ~ d^{sm}_p.\]
\end{lemma}
\proof
If $h(s)$ is a holomorphic function defined in an open neighbourhood
of $s = 0$ then one easily checks using
$\Gamma(s+1) = s \cdot \Gamma(s)$ and $\Gamma(1) = 1$
\begin{eqnarray*}
\left.\frac{d}{ds}\frac{1}{\Gamma(s)} \cdot h(s) \right|_{s=0}
& = &  h(0);
\\
\left.\frac{d}{ds}\frac{1}{\Gamma(s)} \cdot \frac{1}{s}\right|_{s=0}
& = &  - \left.\frac{d\Gamma}{ds}\right|_{s=1}.
\end{eqnarray*}
Notice that the function
$\int_0^1 t^{s-1} \cdot \left( \tr_\Gamma
e^{-t\Delta_p^{\perp}[\widetilde{N}]}-
\sum_{i=0}^m t^{-(m-i)/2}(a_i+b_i)\right) dt$ is holomorphic
for $\Re(s)>-1/2$. Hence the claim follows from the following computation
\begin{eqnarray*}
\left.\frac{d}{ds}\frac{1}{\Gamma(s)}
\int_0^1 t^{s-1} \cdot t^{-(m-i)/2} dt\right|_{s=0}
& = &
\left.\frac{d}{ds}\frac{1}{\Gamma(s)}
\cdot \frac{1}{s - (m-i)/2} \right|_{s=0} .\qed
\end{eqnarray*}

\begin{proposition}\label{dsm} In the situation of Theorem
\ref{Tor_converg} and with Notation \ref{notation for the ends}
we get
\[ \lim_{R \to \infty} d_p^{sm}[\widetilde{M_R}] ~ = ~
d_p^{sm}[\widetilde{M}]. \]
\end{proposition}
\proof
First we choose the fundamental domain ${\cal F}\subset \widetilde{M}$
for the $\Gamma = \pi_1(X)$-action on $\widetilde{M}$
such that
\begin{eqnarray*}
{\cal F}_R  & := & \widetilde{M_R} \cap {\cal F};
\\
\partial{\cal F}_R & := & \partial \widetilde{M_R} \cap {\cal F};
\end{eqnarray*}
are fundamental domains for the induced $\Gamma$-action on
$\widetilde{M_R}$ and $\widetilde{\partial M_R}$ and under the
identifications of Notation \ref{notation for the ends}
we get
\begin{eqnarray}
{\cal F}_R - {\cal F}_S & = & [S,R] \times \coprod_j {\cal G}_j;
\label{eqn 2.30}
\\
\partial {\cal F}_R & = & \{R\} \times \coprod_j {\cal G}_j,
\label{eqn 2.31}
\end{eqnarray}
where ${\cal G}_j \subset \R^{m-1}$ is a fundamental domain for
the universal covering $\R^{m-1} \longrightarrow F_j$ of the flat
closed $m-1$-dimensional manifold sitting at the $j$-th component
of $E_0$.

We will only consider $R>2$. We have to estimate
$\abs{d_p^{sm}[\widetilde{M}]-d_p^{sm}[\widetilde{M_R}]}$.
Recall from Definition \ref{definition of normalized trace}
and Lemma \ref{asymptotic expansion of the heat kernel} that
\[ \begin{split}
\tr_\Gamma e^{-t\Delta^{\perp}_p[\widetilde{M_R}]}&=
 \int_{{\cal F}_R}\Spur_{\R} e^{-t\Delta^{\perp}_p[\widetilde{M_R}]}(x,x)\,dx;\\
 a_i[\widetilde{M_R}] &=\int_{{\cal F}_R}\alpha_i[\widetilde{M_R}](x)\,dx;\\
b_i[\widetilde{M_R}] &=
\int_{\boundary {\cal F}_R}\beta_i[\widetilde{M_R}](x)\,dx
\end{split}\]
and similiarly for $\widetilde{M}$. The functions
$\alpha_i[\widetilde{M_R}](x)$ and $\beta_i[\widetilde{M_R}](x)$ are
determined by the geometry of $\widetilde{M_R}$ in a
neighbourhood of $x$. In particular,
$\alpha_i[\widetilde{M_R}]$ does not depend on $R$, and
coincides with $\alpha_i[\widetilde{M}]$.
Moreover, $\Spur e^{-t\Delta^{\perp}_p[\widetilde{M_R}]}(x,x)-
\sum_{i=0}^m t^{-m/2+i/2}\alpha_i[\widetilde{M_R}](x)$
is $O(t^{1/2})$ uniformly in $x$ for
$x\in \widetilde{M_{R-1}}\subset \widetilde{M_R}$
by Theorem \ref{heat_ker_comp}
applied to $N = \widetilde{M}$, $V = \widetilde{M_R}$ and $K = 1$
since $\alpha_i[\widetilde{M_R}](x) = \alpha_i[\widetilde{M}](x)$.
The function
$\Spur e^{-t\Delta^{\perp}_p[\widetilde{M}]}(x,x)-
\sum_{i=0}^m t^{-m/2+i/2}\alpha_i[\widetilde{M}](x)$
is $O(t^{1/2})$ uniformly in $x$ since
the isometry group of $\widetilde M=\HH^m$ is transitive.
This shows that the following splitting makes sense, i.e.
the integrals do converge and the obvious interchange of integration
are allowed. Namely, we write
$d_p^{sm}[\widetilde{M}]-d_p^{sm}[\widetilde{M_R}]$ as a sum
with the following summands.
\begin{equation*}\begin{split}
s_1  := &
\int_{{\cal F}_{R/2}}\int_0^1
\Bigl(\Spur e^{-t\Delta^{\perp}_p(\widetilde M)}(x,x) -
\Spur e^{-t\Delta^{\perp}_p[\widetilde{M_R}]}(x,x)
\\
& \hspace{8mm}-  \sum_i t^{-m/2+i/2}
 (\underbrace{\alpha_i[\widetilde{M}](x)-
\alpha_i[\widetilde{M_R}](x)}_{=0})\Bigr)\frac{dt}{t}\,dx;
\\
s_2   := &  \int_{{\cal F}_{R-1}-{\cal F}_{R/2}}
 \int_0^1\Bigl(\Spur e^{-t\Delta^{\perp}_p(\widetilde M)}(x,x)
 -\Spur e^{-t\Delta^{\perp}_p[\widetilde{M_R}]}(x,x)
\\
 & \hspace{8mm}  -\sum_i
t^{-m/2+i/2}
 (\underbrace{\alpha_i[\widetilde{M}](x)-
\alpha_i[\widetilde{M_R}](x)}_{=0})\Bigr)\frac{dt}{t}\,dx;
\\
s_3  := &
\int_{{\cal F}-{\cal F}_{R-1}}\int_0^1
 \left(\Spur e^{-t\Delta^{\perp}_p(\widetilde M)}(x,x)-
 \sum_i t^{-m/2+i/2}\alpha_i[\widetilde{M}](x)\right)\frac{dt}{t}\,dx;
\\
s_4  :=  & \int_0^1\left(
\int_{{\cal F}_R-{\cal F}_{R-1}}\Spur
e^{-t\Delta^{\perp}_p[\widetilde{M_R}]}(x,x)
\,dx\right.
\\
 & \hspace{8mm} \left.-  \sum_i t^{-m/2+i/2}\left(
 \int_{{\cal F}_R-{\cal F}_{R-1}}\alpha_i[\widetilde{M_R}](x)\,dx +
\int_{\boundary {\cal F}_R}\beta_i[\widetilde{M_R}](x')\,dx'\right)
\right)\frac{dt}{t};
\\
s_5  := & 
\sum_{i=0}^m c(i,m)\int_{{\cal F}_R}
 \underbrace{\alpha_i[\widetilde{M}](x)-
\alpha_i[\widetilde{M_R}](x)}_{=0}\,dx ~ = ~ 0;
\\
s_6  := &
\sum_{i=0}^m c(i,m)
\int_{{\cal F}-{\cal F}_R}\alpha_i[\widetilde{M}](x)\,dx;
\\
s_7  := & 
\sum_{i=0}^m c(i,m)\int_{\boundary {\cal
F}_R}\beta_i[\widetilde{M_R}](x')\,dx'.
\end{split}\end{equation*}
We study each of these summands individually.
For $s_1$ and $s_2$ we use Theorem \ref{kernComp}
applied to $N = \widetilde{M}$, $V = \widetilde{M_R}$ and $K = 1$.
Note that $d(\widetilde M_a, \widetilde M-\widetilde M_b)=b-a$
for $b/2\le a\le b$. This implies for appropriate constants
$C_1$ and $C_2$ independent of $R$:
\[\begin{split}
  \abs{s_1} &\le \vol(M_{R/2})
  \int_0^1 C_1 e^{-R^2/4C_2 t}\frac{dt}{t}\le
  4\vol (M) C_1 C_2 R^{-2} e^{-R^2/4C_2}; \\
  \abs{s_2} &\le \vol(M_{R-1}-M_{R/2})
  \int_0^1 C_1 e^{-1/C_2 t}\frac{dt}{t} \le
  \vol(M-M_{R/2}) C_1 C_2 =\vol(E_{R/2}) C_1 C_2.
\end{split}\]
Therefore, both terms tend to zero for $R\to\infty$.
For $s_3$ and $s_6$ observe that
$\widetilde M=\HH^m$ has transitive isometry group. It
follows that
\[ I_3:=\int_0^1\left(\Spur e^{-t\Delta^{\perp}_p(\widetilde M)}(x,x)-
 \sum_i t^{-m/2+i/2}\alpha_i[\widetilde{M}](x)\right)\frac{dt}{t} \]
is a constant independent of $x$ and the same holds for
\[ I_6:=\sum_{i=0}^m c(i,m)\alpha_i[\widetilde{M}](x). \]
 Therefore
\[\begin{array}{lclclcl}
\abs{s_3}  & \le &
\abs{I_3}\cdot \vol(M-M_{R-1}) & = & \abs{I_3}\cdot \vol(E_{R-1})
& \xrightarrow{R\to\infty} &0; \\
 \abs{s_6} & \le & \abs{I_6}\cdot \vol(M-M_R) &
 = & \abs{I_6}\cdot \vol(E_R ) & \xrightarrow{R\to\infty} & 0.
\end{array} \]
For arbitrary $R,S\ge 0$ and
 $x'\in\boundary \widetilde M_R$ and $y'\in\boundary\widetilde M_S$ we find
neighbourhoods which are isometric. It follows that
\[ I_7:=\sum_{i=0}^m c(i,m)\beta_i[\widetilde{M_R}](x') \]
does not depend on $x'$ neither on $R$. Therefore
\[ \abs{s_7}\le \abs{I_7}
\vol(\boundary M_R)\xrightarrow{R\to\infty} 0. \]

Since $s_5$ is zero it remains to treat $s_4$.
We will treat only the case where
$M_R - M_{R-2}$ has only one component,
otherwise one applies the following argument to each component
separately. Define
$$\HH^m_R ~ := ~ (-\infty,R] \times \R^{m-1}$$
with the warped product metric as
in Notation \ref{notation for the ends}.
We split $s_4$ into three summands
\begin{eqnarray*}
s_{41}& := &\int_0^1 \int_{{\cal F}_R-{\cal F}_{R-1}}
\left(\Spur e^{-t\Delta^{\perp}_p[\widetilde{M_R}]}(x,x)
  -\Spur e^{-t\Delta^{\perp}_p[\widetilde{M_R} - \widetilde{M_{R-2}}]}(x,x)
  \right)\,dx \frac{dt}{t};
\\
s_{42}& := &\int_0^1 \int_{{\cal F}_R-{\cal F}_{R-1}}
\left(\Spur e^{-t\Delta^{\perp}_p[\widetilde{M_R} - \widetilde{M_{R-2}}]}(x,x)
  -\Spur e^{-t\Delta^{\perp}_p[\HH^m_R ]}(x,x)
  \right)\,dx \frac{dt}{t};
\\
s_{43}& :=  & \int_0^1 \int_{{\cal F}_R-{\cal F}_{R-1}}
\Spur e^{-t\Delta^{\perp}_p[\HH^m_R]}(x,x) \,dx
\\ & & \hspace{10mm}
- \left(\sum_i t^{-m/2+i/2}\int_{{\cal F}_R-{\cal F}_{R-1}}
\alpha_i[\widetilde{M_R}](x)\,dx +
\int_{\boundary {\cal
F}_R}\beta_i[\widetilde{M_R}](x')\,dx'\right)\frac{dt}{t}.
\end{eqnarray*}
For $s_{41}$ we want to apply Theorem
\ref{kernComp} for $V = \widetilde{M_R} - \widetilde{M_{R-2}}$,
$N = \widetilde{M_R}$ and $K = 1$. Since $\widetilde{M_R}$ has constant
sectional curvature $-1$ for all $R$ and a neighbourhood of $\widetilde{M_R}$
is isometric to a neighbourhood of $\widetilde{M_0}$, the constants appearing
in Theorem  \ref{kernComp} can be chosen independently of $R$.

If $x \in \widetilde M_R-\widetilde M_{R-1}$
then $d_{\widetilde{M_R} -\widetilde{M_{R-2}}}(x) \ge 1$ and
Theorem \ref{kernComp} yields
\[
\abs{e^{-t\Delta^{\perp}_p[\widetilde{M_R}]}(x,x) -
e^{-t\Delta^{\perp}_p[\widetilde{M_R} - \widetilde{M_{R-2}}]}(x,x)}
\le D_1 \cdot e^{-D_2/t}\]
for constants $D_1$ and $D_2$ independent of $x$, $t$ and $R$.
Since the volume of $\widetilde{M_R} -\widetilde{M_{R-1}}$
tends to zero if $R$ goes to $\infty$ the same is true for
$s_{41}$. The same argument when replacing
$\widetilde{M_R}$ by $\HH_R^m$ yields
that $s_{42}$ tends to zero if $R$ goes to $\infty$.

Recall that $\alpha[\widetilde{M_R}](x,x)$ and
$\beta[\widetilde{M_R}](x,x)$ are determined by the geometry
of $\widetilde{M_R}$ in a neighborhood of $x$.
Hence we conclude
\begin{eqnarray*}
s_{43}& :=  & \int_0^1 \int_{{\cal F}_R-{\cal F}_{R-1}}
\Spur e^{-t\Delta^{\perp}_p[\HH^m_R]}(x,x) \,dx
\\ & & \hspace{10mm}
- \left(\sum_i t^{-m/2+i/2}\int_{{\cal F}_R-{\cal F}_{R-1}}
\alpha_i[\HH_R^m](x)\,dx +
\int_{\boundary {\cal F}_R}\beta_i[\HH_R^m](x')\,dx'\right)\frac{dt}{t}.
\end{eqnarray*}
In the sequel we use the identifications
\eqref{eqn 2.30} and \eqref{eqn 2.31}.
Notice that each isometry $i: \R^{m-1} \longrightarrow \R^{m-1}$
induces an isometry $\id \times i: \HH^m_R \longrightarrow \HH^m_R$
and that $\HH^m_R$ is isometric to $\HH_m^0$ by the same argument
as in the proof of Lemma \ref{isometries}. This implies
that
\begin{eqnarray*}
\beta_i[\HH_R^m](R,y) & = & f_i;
\\
\alpha_i[\HH_R^m](u,y) & = & g_i(u + 1 - R);
\\
\Spur e^{-t\Delta^{\perp}_p[\HH^m_R]}((u,y),(u,y)) & = & h(t,u + 1- R);
\end{eqnarray*}
holds for all $u \in (-\infty,R]$,
$y \in \R^{m-1}$ and
an appropriate number $f_i$, appropriate functions $g_i(u)$ and an
appropriate function $h(t,u)$,
which are all independent of $y$ or $R$.
Since the volume of $\{u\} \times {\cal G}$ in $\{u\} \times \R^{m-1}$
is $e^{-(m-1)u}$-times the volume of ${\cal G} \subset \R^{m-1}$, we get
\begin{eqnarray*}
s_{43}& = & e^{-(m-1)R} \cdot \vol({\cal G})\cdot
 \int_0^1 \int_0^1\left(  h(t,u) e^{-2u+2}\,du
- \sum_i t^{-\frac{m+i}{2}}\cdot  \left(g_i(u) e^{-2u+2} -
f_i \right)\right) du\frac{dt}{t}.
\end{eqnarray*}
Hence $s_{43}$ tends to zero if $R$ goes to $\infty$.
This finishes the proof
of Proposition \ref{dsm}. \qed



\typeout{--------------------- section 3-----------------------}

\section{The large t summand of the torsion
corresponds to small eigenvalues}
\label{The large t summand of the torsion corresponds to
small eigenvalues}

Here, we will show that it sufficies to
control uniformly the small eigenvalues,
to get convergence of
$\int_1^\infty t^{-1}
\tr_\Gamma e^{-t\Delta^{\perp}_p[\widetilde{M_R}]}\,dt$.
For this section let $p$
be a fixed non-negative integer. We
start with some notation.

\begin{definition}
Let $N$ be a compact Riemannian manifold
possibly with boundary.
Let $E^p_{\lambda}[\widetilde{N}]$ be
the right-continuous spectral family
of $\Delta^{\perp}_p[\widetilde{N}]$, i.e.
$\Delta_p[\widetilde{N}]$ restricted to the
complement of its kernel. Define
\begin{eqnarray}
F(\Delta^{\perp}_p[\widetilde{N}],\lambda) & = &
\tr_\Gamma E^p_{\lambda}[\widetilde{N}].
\label{eqn 3.1}
\end{eqnarray}
Abbreviate
$F(\lambda) = F(\Delta^{\perp}_p[\widetilde{N}],\lambda)$
and $E_{\lambda} = E^p_\lambda[\widetilde{N}]$. \qed
\end{definition}
Definition \ref{eqn 3.1} is of course consistent with the definition
(see Definition \ref{definition of spectral density functions}
and Lemma \ref{compute})
of spectral density function we will use later.

\begin{lemma}
For every $\lambda<\infty$ we have $F(\lambda)<\infty$.
\end{lemma}
\begin{proof}
  Remember that $e^{-\Delta_p[\tilde N]}$ and therefore also
$e^{-\Delta_p^\perp[\tilde N]}$
 are of $\Gamma$-trace class. It follows that the
operator
$\chi_{[0,\lambda]}(\Delta_p^\perp)e^{-\Delta_p^\perp}$ has the same property,
since the characteristic function $\chi_{[0,\lambda]}(t)$
of the interval $[0,\lambda]$ is clearly bounded and the trace class operators
form an ideal in the algebra of bounded operators. Then
$E(\lambda)=\chi_{[0,\lambda]}(\Delta_p^\perp)e^{\Delta_p^\perp}
\chi_{[0,\lambda]}(\Delta_p^\perp)e^{-\Delta_p^\perp}$ is also of
$\Gamma$-trace class because $\chi_{[0,\lambda]}(t)e^t$ is a bounded function,
too. This concludes the proof.
\end{proof}

Now observe that for $t\ge 1$
\begin{eqnarray}
t^{-1}\tr_\Gamma e^{-t\Delta^{\perp}_p[\widetilde{N}]}
& = &
t^{-1}\int_0^\epsilon e^{-t\lambda}
\,dF(\lambda) + t^{-1}\tr_\Gamma\int_\epsilon^\infty e^{-t\lambda}
\,d E_{\lambda}
\nonumber
\\
& \stackrel{F(0)=0}{=} & -t^{-1}\int_0^\epsilon(-t)e^{-t\lambda}
F(\lambda)\,
d\lambda +t^{-1}e^{-t\epsilon}F(\epsilon)
\nonumber
\\
& & \hspace{30mm}
+ t^{-1}
e^{-t\epsilon}\tr_\Gamma\int_\epsilon^\infty
e^{-t(\lambda-\epsilon)}\,d E_{\lambda}
\nonumber
\\
& \stackrel{t\ge 1}{\le} &
\int_0^\epsilon e^{-t\lambda}
F(\lambda)
\,d\lambda +
\frac{e^{-t\epsilon}}{t}F(\epsilon) +\frac{e^{-t\epsilon}}{t}
\tr_\Gamma\int_\epsilon^\infty e^{-(\lambda-\epsilon)}d E_{\lambda}
\nonumber
\\
&\le &
\int_0^\epsilon e^{-t\lambda}F(\lambda)\,d\lambda +
\frac{e^{-t\epsilon}}{t}F(\epsilon) +
\frac{e^{-t\epsilon}}{t}e^\epsilon
\tr_\Gamma\int_0^\infty e^{-\lambda} d E_{\lambda}
\nonumber
\\
& = &
\int_0^\epsilon e^{-t\lambda}F(\lambda)\,d\lambda +
\frac{e^{-t\epsilon}}{t}F(\epsilon) +\frac{e^{-t\epsilon}}{t}e^\epsilon
\tr_\Gamma e^{-\Delta^{\perp}_p[\widetilde{N}]}.
\label{ew}
\end{eqnarray}

Therefore the next step in the proof of
Theorem \ref{Tor_converg} is
 the following:
\begin{proposition}\label{const_t_conv}
For every $t\ge 1$ we have
\[ \lim_{R \to \infty}~ \tr_\Gamma e^{-t\Delta^{\perp}_p[\widetilde{M_R}]}
~ = ~
 \tr_\Gamma e^{-t\Delta^{\perp}_p[\widetilde{M}]}. \]
\end{proposition}
\proof
By  Theorem \ref{ker_est} and the local isometries of
Lemma \ref{isometries} there is $c=c(1)$
independent of $R$ so that
\[ \abs{e^{-t\Delta_p[\widetilde{M_R}]}(x,x)}
\le c\qquad\forall t\ge 1\;
\forall x\in\widetilde{M_R}. \]
Since $\widetilde{M}=\HH^m$ is homogeneous we can choose $c$ so that this
inequality also holds for $M$ in place of $M_R$.
Notice that $m$ is odd
by assumption. Hence
$\Delta^{\perp}_p[\widetilde{M}] =\Delta_p[\widetilde{M}]$
by \cite{Dodziuk (1979)}. We conclude
$\Delta^{\perp}_p[\widetilde{M_R}] =\Delta_p[\widetilde{M_R}]$
from Corollary \ref{L2_acyclic} or
\cite[Theorem 1.1]{Cheeger-Gromov (1985)}.
Given $t \ge 1$, we get for $R\ge 2$
\begin{eqnarray*}
\lefteqn{\abs{\tr_\Gamma e^{-t\Delta^{\perp}_p[\widetilde{M_R}]}-
\tr_\Gamma e^{-t\Delta^{\perp}_p[\widetilde{M}]}}}
\\
&  = &
\abs{\int_{{\cal F}_R}\Spur e^{-\Delta_p[\widetilde{M_R}]}(x,x)\,dx
-\int_{{\cal F}}\Spur e^{-\Delta_p(\widetilde{M}]}(x,x)\,dx }
\\
&\le & \int_{{\cal F}_{R/2}}
\abs{\Spur e^{-t\Delta_p[\widetilde{M_R}](x,x)} -
\Spur e^{-t\Delta_p(\widetilde{M}]}(x,x)}
\,dx
\\
& & \qquad + \int_{{\cal F}_R-{\cal F}_{R/2}}
\Spur e^{-t\Delta_p[\widetilde{M_R}]}(x,x)\,dx +
\int_{{\cal F} - {\cal F}_{R/2}}
\Spur e^{-t\Delta_p[\widetilde{M}])}(x,x)\,dx
\\
& \le & \vol(M_{R/2})C_1e^{-R^2/4C_2t}
+ c\vol(M_R-M_{R/2}) + c\vol(E_{R/2})
\\
&\le & \vol(M)C_1e^{-R^2/4C_2t}  + 2c\vol(E_{R/2})
\xrightarrow{R\to\infty} 0.
\end{eqnarray*}
The existence of $C_1, C_2 > 0$ follows from Theorem
\ref{heat_ker_comp} applied to $N = \widetilde{M}$,
$V = \widetilde{M_R}$, $K = R/2$. Note that
$d(\tilde M_{R/2},\tilde M-\tilde M_R)= R/2$.
This finishes the proof of Proposition \ref{const_t_conv}.\qed

\begin{corollary}\label{goal}
If we find $\epsilon>0$ and a function
$G(\lambda)$ so that for every $R\ge0$
\[ F(\Delta^{\perp}_p[\widetilde{M_R}],\lambda)\le G(\lambda)
  \quad\forall \lambda\le \epsilon \]
with
\[ \int_1^\infty
 \left(\int_0^\epsilon e^{-t\lambda}G(\lambda)\,d\lambda\right)\,
dt <\infty ,\]
then the large time summand in the
 analytic $L^2$-torsion of $M_R$ converges to
the corresponding summand for $M$
\[ \int_1^\infty t^{-1}\tr_\Gamma e^{-t\Delta^{\perp}(\widetilde{M_R})}\,dt
\xrightarrow{R\to\infty}
\int_1^\infty t^{-1}\tr_\Gamma e^{-t\Delta^{\perp}(\widetilde{M})}\,dt. \]
\end{corollary}
\proof
By Proposition \ref{const_t_conv} we
have pointwise convergence of the integrand.
We want to apply the Theorem of Dominated
Convergence. Inequality \eqref{ew} shows
that the assumption just guaranties the
existence of a dominating integrable
function because $\tr_\Gamma e^{-\Delta^{\perp}(\widetilde{M_R})}$
is bounded by Proposition \ref{const_t_conv}
and $\int_1^\infty e^{-t\epsilon}/t\,dt <\infty$.
\qed



\typeout{--------------------- section 4 -----------------------}

\section{Spectral density functions}
\label{Spectral density functions}

In this section we have to go through some of the proofs
of \cite[Section 1 and Section 2]{Lott-Lueck (1995)} since
we need the results there in a more precise form later.
For this section, let $\Neumann{A}$
be a von Neumann algebra with finite trace $\tr$.
Our main example will be
$\Neumann{A}=\Neumann{N}(\Gamma)$ the von Neumann
algebra of a discrete group. In this section,
all morphisms will be
Hilbert-$\Neumann{A}$-module morphisms, i.e. bounded
$\Neumann{A}$-equivariant
operators, unless explicitly specified differently.

\subsection{Spectral density functions}

\begin{definition} \label{definition of spectral density functions}
Suppose $f:U\to V$ is a
Hilbert-$\Neumann{A}$-module morphism. Define its
{\em spectral density function}
\[ F(f,\lambda):= \tr E_{\lambda^2}(f^*f)\qquad\lambda\ge 0. \]
Here $E_\lambda(f^*f)$ is the right-continuous spectral
family of the positive self adjoint
operator $f^*f$. Note that $F(f,\lambda)$ is monoton increasing.
We say $f$ is {\em left Fredholm} if $F(f,\lambda)<\infty$
for some $\lambda>0$. If $F(f,0)<\infty$ set
\[\overline{F}(f,\lambda):=F(f,\lambda)-F(f,0).\]
\end{definition}

\begin{lemma}\label{basic_F}
Let $f:U\to V$ and $g:V\to W$ be given.
Let $i:V\to V'$ be injective with
closed range and $p:U\to U'$ surjective
with $\dim_{\Neumann{A}}(\ker(p))<
\infty$. Then
\begin{enumerate}
\item $F(f,\lambda)\le F(gf,\norm{g}\lambda)\qquad\forall\lambda$;

\item $F(g,\lambda)\le F(gf,\norm{f}\lambda)$, if $f$ is left
Fredholm with dense image;

\item $F(gf,\lambda)\le
F(g,\lambda^{1-r})+F(f,\lambda^r)$ $\forall 0<r<1$;

\item $ F(if,\lambda)\le F(f,\norm{i^{-1}}\lambda)$,
where $i^{-1}:i(V)\to V $
is bounded by the Open Mapping Theorem;

\item $F(f,\lambda)\le F(fp,\norm{p}\lambda)$;

\item $F(f^*f,\sqrt{\lambda})=F(f,\lambda)$.
\end{enumerate}
\end{lemma}
\proof 1.) - 3.) are proven in \cite[Lemma 1.6]{Lott-Lueck (1995)}
and imply 4.) and 5.). 6.) is a direct consequece of the definition. \qed

\begin{lemma}\label{basic_barF}
Adopt the situation of Lemma \ref{basic_F}
and suppose that the kernels of
all morphisms in question are finite
$\Neumann{A}$-dimensional. Then
\begin{enumerate}
\item  $\overline{F}(f,\lambda)\le \overline{F}(gf,\norm{g}\lambda)$,
if $\ker g\cap\im f=\{0\}$;

\item $\overline{F}(g,\lambda)\le \overline{F}(gf,\norm{f}\lambda)$,
if $f$ is left Fredholm with dense image;

\item $\overline{F}(gf,\lambda)\le \overline{F}(g,\lambda^{1-r})+
\overline{F}(f,\lambda^r)$ for all $0<r<1$,
if $\ker g \subset\overline{\im f}$;

\item $ \overline{F}(if,\lambda)\le \overline{F}(f,\norm{i^{-1}}\lambda)$;

\item  $\overline{F}(fp,\lambda)\le
\overline{F}(f,\norm{p^{-1}}\lambda)$, for
$p^{-1}:U\to (\ker p)^\perp$;

\item $\overline{F}(f,\lambda)\le
\overline{F}(fp,\norm{p}\lambda)+\dim_\Neumann{A}\ker p$.
\end{enumerate}
\end{lemma}
\proof This follows from Lemma \ref{basic_F}. In assertion
2.) use \cite[Lemma 1.4]{Lott-Lueck (1995)} to conclude
$F(gf,0) \le F(g,0)$.
In assertion 3.) use the easy argument in the proof of
\cite[Lemma 1.11.3]{Lott-Lueck (1995)} to conclude
$F(gf,0) = F(g,0) + F(f,0)$. \qed

\begin{lemma} \label{F and adjoint}
If $\dim_{\Neumann{A}}\ker\phi<\infty$ and
$\dim_{\Neumann{A}}\ker\phi^*<\infty$, then
\[ \overline{F}(\phi,\lambda) =\overline{F}(\phi^*,\lambda) \]
This also holds if $\phi$ is an
unbounded $\Neumann{A}$-operator.
\end{lemma}
\proof
The proof in \cite[Lemma 1.12.6]{Lott-Lueck (1995)}
literally holds for unbounded operators, too. \qed

\begin{lemma} \label{F and block matrices}
Let $\phi:U_1\to V_1$, $\gamma:U_2\to V_1$
and $\xi:U_2\to V_2$ be morphisms
of Hilbert-$\Neumann{A}$-modules. Then
\begin{enumerate}
\item $ F\left(\left(\begin{smallmatrix}
\phi & 0\\ 0& \xi
\end{smallmatrix}\right),\lambda\right) =
F(\phi,\lambda)+ F(\xi,\lambda)$;
\item Suppose $\phi$ is invertible. Then
\[ F\left(\begin{pmatrix}
    \phi & \gamma\\ 0& \xi
  \end{pmatrix},\lambda\right)
\le F(\phi,(4 +2\norm{\gamma}\norm{\phi^{-1}})\lambda) +
 F(\xi,(4 +2\norm{\gamma}\norm{\phi^{-1}})\lambda);\]
\item $ F\left(\begin{pmatrix}\phi & \gamma\\ 0& \xi
\end{pmatrix},\lambda\right)
\le F(\phi,\lambda^{r})+
F(\xi,(4+2\norm{\gamma})\lambda^{1-r})$
holds for $0<r<1$, provided
that $\lambda < (4+2\norm{\gamma})^{1/(r-1)}$ is true;

\item
$F(\phi,\lambda)\le  F\left(\begin{pmatrix}
\phi & \gamma\\ 0& \xi
\end{pmatrix},2(1+\norm{\gamma}+\norm{\xi})\lambda\right);$

\item If $\phi$ has dense image and $\phi$ or
$\left(\begin{smallmatrix}\phi&\gamma\\0& \xi \end{smallmatrix}
\right)$ are left Fredholm
then for $\lambda < 1$ we have
\[ F(\xi,\lambda)\le  F\left(\begin{pmatrix}
\phi & \gamma\\ 0& \xi
\end{pmatrix},2(1+\norm{\gamma}+\norm{\phi})\lambda\right). \]
\end{enumerate}
\end{lemma}
\proof
We will use the elementary fact
\[ \left\|\begin{pmatrix}
\phi & \gamma\\ 0& \xi
\end{pmatrix}\right\| \le 2(\norm{\phi}+\norm{\xi}+\norm{\gamma}) \]
and the decompositions
\begin{eqnarray}
\begin{pmatrix}  \phi & 0 \\ 0& \xi   \end{pmatrix}
& = &
\begin{pmatrix} \phi&\gamma\\  0&\xi \end{pmatrix}
\begin{pmatrix} 1& - \phi^{-1}\gamma \\ 0&1  \end{pmatrix};
\label{eqn 4.1}
\\
\begin{pmatrix} \phi &\gamma\\ 0&\xi \end{pmatrix}
&  = &
\begin{pmatrix} 1 &\gamma\\ 0&\xi \end{pmatrix}
\begin{pmatrix} \phi& 0\\ 0&1  \end{pmatrix};
\label{eqn 4.2}
\\
\begin{pmatrix} \phi &\gamma\\ 0&\xi \end{pmatrix}
&  = &
\begin{pmatrix} 1&0\\ 0&\xi  \end{pmatrix}
\begin{pmatrix} \phi&\gamma\\ 0&1 \end{pmatrix}.
\label{eqn 4.3}
\end{eqnarray}
1.) is obvious.\par\noindent
2.) Apply Lemma \ref{basic_F}.5 to \eqref{eqn 4.1} and use
assertion 1.)\par\noindent
3.) Apply Lemma \ref{basic_F}.3 to \eqref{eqn 4.2}
and use assertions 1.) and 2.).\par\noindent
4.) Apply Lemma \ref{basic_F}.1 to \eqref{eqn 4.2}. \par\noindent
5.) If
$\left(\begin{smallmatrix}\phi&\gamma\\0& \xi \end{smallmatrix}\right)$
is left Fredholm, then $\phi$ is left Fredholm
by \cite[Lemma 1.11.1]{Lott-Lueck (1995)}. Suppose
that $\phi$ has dense image and is left Fredholm. Then
$\left(\begin{smallmatrix}\phi&\gamma\\0& 1 \end{smallmatrix}\right)$
has dense image and is left Fredholm by
\cite[Lemma 1.12.3]{Lott-Lueck (1995)}.
Apply Lemma \ref{basic_barF}.2 to \eqref{eqn 4.3}
and use assertion 1.). \qed

\begin{lemma}\label{matrix_spec_dens}
Adopt the situation of Lemma \ref{F and block matrices}.
Suppose all relevant kernels have
finite $\Neumann{A}$-dimension. Then
\begin{enumerate}
\item $ \overline{F}\left(\begin{pmatrix}
\phi & 0\\ 0& \xi
\end{pmatrix}\,,\lambda\right) = \overline{F}(\phi,\lambda)+
\overline{F}(\xi,\lambda)$;

\item Suppose $\phi$ is invertible. Then
\[ \overline{F}\left(\begin{pmatrix}
\phi & \gamma\\ 0& \xi
\end{pmatrix},\lambda\right) \le
\overline{F}(\phi,(4 +2\norm{\gamma}\norm{\phi^{-1}})\lambda)+
\overline{F}(\xi,(4 +2\norm{\gamma}\norm{\phi^{-1}})\lambda);\]

\item Suppose $\xi$ is injective or $\phi$ has dense image and is left
Fredholm. Then for $0 < r < 1$ and $\lambda < (4+2\norm{\gamma})^{1/(r-1)}$
we have
\[ \overline{F}\left(\begin{pmatrix}
\phi & \gamma\\ 0& \xi
\end{pmatrix},\lambda\right) \le \overline{F}(\phi,\lambda^{r})+
\overline{F}(\xi,(4 + 2\norm{\gamma})\lambda^{1-r}).\]

\item If $\xi$ is injective then
\[ \overline{F}(\phi,\lambda)\le  \overline{F}\left(\begin{pmatrix}
\phi & \gamma\\ 0& \xi
\end{pmatrix},2(1+\norm{\gamma}+\norm{\xi})\lambda\right); \]

\item If $\phi$ has dense image and $\phi$ or
$\left(\begin{smallmatrix}\phi&\gamma\\0& \xi
\end{smallmatrix}\right)$ are left Fredholm,
then
\[ \overline{F}(\xi,\lambda)\le \overline{F}\left(\begin{pmatrix}
\phi & \gamma\\ 0& \xi
\end{pmatrix},2(1+\norm{\gamma}+\norm{\phi})\lambda\right) +
 \dim_{\Neumann{A}}\ker\phi. \]
\end{enumerate}
\end{lemma}
\proof
This follows from Lemma \ref{F and block matrices}
as Lemma \ref{basic_barF} follows from Lemma \ref{basic_F}. \qed

Next we extend the notion of $F(f,\lambda)$ and $\bar F(f,\lambda)$ from
Definition \ref{definition of spectral density functions} for morphisms
to chain complexes.
\begin{definition}
Let $0\to C^0\xrightarrow{c^0} C^1 \xrightarrow{c^1}\cdots $
be a cochain complex of Hilbert-$\Neumann{A}$-modules.
Define the {\em spectral density function}
\[ F^p(C^*,\lambda):= F(c^p|: \im(c^{p-1})^\perp \to C^{p+1},\lambda). \]
If $F^p(C^*,0)<\infty$ set
\[ \overline{F}^p(C^*,\lambda):= F^p(C^*,\lambda)- F^p(C^*,0). \]
We call $C^*$ {\em Fredholm} at $p$ if $F^p(C^*,\lambda)<\infty$ for
some $\lambda>0$. \qed
\end{definition}

\subsection{Short exact sequences of cochain complexes}

In this subsection we will express for a short exact sequence of Hilbert
$\Neumann{A}$-cochain complexes $0\to C^*\to D^*\to E^* \to 0$
the spectral density function of $D^*$ in
terms of the spectral density functions
of $C^*$, $E^*$ and the long weakly exact homology sequence.

\begin{theorem}\label{short_exact}
Let $0\to C^*\xrightarrow{j} D^*\xrightarrow{q} E^*\to 0$
 be an exact sequence of Hilbert cochain
complexes as above. Suppose $C^*$ and $E^*$ are Fredholm at $p$.
Then $D^*$ is Fredholm at $p$, $\delta_p$ is left Fredholm and
\[ \overline{F}_p(D^*,\lambda)\le \overline{F}_p(E^*, c_E\cdot \lambda^{1/2}) +
  \overline{F}(\delta^p, c_\delta\cdot\lambda^{1/4}) +
  \overline{F}_p(C^*, c_C\cdot \lambda^{1/4})
  \qquad \text{for }0\le\lambda<c_1. \]
Here $\delta^p: H^p(E^*)\to H^{p+1}(C^*)$
is the connecting homomorphism in the
long weakly exact $L^2$-cohomology sequence
(\cite[Theorem 2.1]{Cheeger-Gromov (1985)},
\cite[Theorem 2.2]{Lott-Lueck (1995)})
and $c_E$, $c_C$, $c_\delta$, $c_1$ are
explicitely determined in terms of the
norms of $d^p$, $j^p$, $q^p$ and their inverses. We use the convention that
$H^p(E^*)$ and $H^p(C^*)$ are subquotients of the corresponding cochain
complexes with the induced norm. Explicitly:
\begin{eqnarray*}
c_E & = &
(4+2\norm{d^p}) \norm{q_{p+1}}\norm{q_p^{-1}};
\\
c_C & = & \norm{j_{p+1}^{-1}}^{1/2}\norm{j_p};
\\
c_\delta & = & \norm{j_{p+1}^{-1}}^{1/2}
(4 + 2\norm{j_{p+1}^{-1}}\norm{d^p}) \norm{q_p^{-1}};
\\
c_1 & = & \min\{
(4 + 2\norm{d^p})^{-1/2},(4 + 2\norm{j_{p+1}^{-1}}\norm{d^p})^{-1/2}
\}.
\end{eqnarray*}
Here, the inverse of an operator with
closed image always means the obvious operator
from the image to the orthogonal complement of the kernel.
\end{theorem}
\proof
For the proof, we repeat the proof of \cite[Theorem 2.3]{Lott-Lueck (1995)}
(where chain and not cochain complexes are treated)
and take care not only of the
Novikov-Shubin invariants but of all of the
spectral density functions.
Lott-L\"uck  \cite[page 28]{Lott-Lueck (1995)}
construct a commutative diagram
\[\begin{CD} 0 @>>> \ker\overline{q_p} @>i>>
D^p/\overline{\im(d_{p-1})} @>{\overline{q_p}}>>
  E^p/\ker(e^p) @>>> 0\\
&& @VV{\partial_p}V @VV\overline{d^p}V @VV{\overline{e^p}}V\\
0 @>>> C^{p+1}@>{j_{p+1}}>> D^{p+1} @>{q_{p+1}}>> E^{p+1} @>>>0
\end{CD}\]
and show that $\partial_p$ is Fredholm.
The diagram yields a splitting
\begin{equation*}
 D^p/\overline{\im(d_{p-1})}= \ker \overline{q_p}\oplus \ker\overline{q_p}^\perp
\xrightarrow{\overline{d^p}=\left(
\begin{smallmatrix} j_{p+1}\partial_p& \gamma\\
0 & q_{p+1}^{-1}\overline{e^p}\overline{q_p}\end{smallmatrix}\right)}
j(C^{p+1})\oplus\ker q_{p+1}^\perp =D^{p+1}.
\end{equation*}
Because $\norm{\gamma}\le\norm{\overline{d^p}}= \norm{d^p}$
and $q_{p+1}^{-1}\overline{e^p}\overline{q_p}$ is injective,
Lemma \ref{matrix_spec_dens}.3 implies
\begin{eqnarray*}
\overline{F}_p(D^*,\lambda) & = &\overline{F}(\overline{d^p},\lambda)
\\
& \le & \overline{F}(j_{p+1}\partial_p,\lambda^{1/2})
+
\overline{F}(q_{p+1}^{-1}\overline{e^p}\overline{q_p},(4+2\norm{d^p})\lambda^{1/2}
)
\qquad\forall \lambda< (4 + 2\norm{d^p})^{-2}.
\end{eqnarray*}
Since $\norm{\overline{q_p}^{-1}}\le\norm{q_p^{-1}}$ we conclude
from Lemma \ref{basic_barF} 4.) and 5.)
 that
for all $\lambda < (4 + 2\norm{d^p})^{-2}$
\begin{eqnarray}
\overline{F}_p(D^*,\lambda)
&   \le &
\overline{F}(\partial_p,\norm{j_{p+1}^{-1}}\lambda^{1/2}) +
 \overline{F}(\overline{e^p},(4+2\norm{d^p})
 \norm{q_{p+1}}\norm{\overline{q_p}^{-1}}\lambda^{1/2})
\nonumber \\
& \le & \overline{F}(\partial_p,\norm{j_{p+1}^{-1}}\lambda^{1/2}) +
 \overline{F}_p(E^*,(4+2\norm{d^p})
 \norm{q_{p+1}}\norm{q_p^{-1}}\lambda^{1/2}).
 \label{eqn 4.10}
\end{eqnarray}

Now we have to examine $\partial_p$ further. Its range actually
lies  in $\ker(c_{p+1})$ \cite[page 28]{Lott-Lueck (1995)}
and Lemma \ref{basic_barF}.4 applied to
$\ker(\overline{q_p})\xrightarrow{\partial_p} \ker(c_{p+1})
\stackrel{i}{\hookrightarrow} C^{p+1}$ implies
\begin{eqnarray}
\overline{F}(\partial_p,\lambda) & = &
\overline{F}(i\circ \partial_p,\lambda). \label{eqn 4.11}
\end{eqnarray}
Lott-L\"uck \cite[page 29]{Lott-Lueck (1995)}
construct the following commutative diagram with exact rows
\[\begin{CD}
 0 @>>> \ker(\widehat{q_p}) @>{i_1}>> \ker(\overline{q_p})
 @>\widehat{q_p}>> H_p(E) @>>> 0\\
 &&    @VV{\overline{\partial_p}}V @VV{\partial_p}V @VV{\delta_p}V\\
0 @>>> \overline{\im(c^p)} @>{i_2}>> \ker(c_{p+1}) @>{pr}>> H_{p+1}(C) @>>> 0
\end{CD}\]
and prove that the induced
operator $\overline{\partial_p}$ is Fredholm and
has dense image. We get the splitting
\begin{multline*}
 \ker(\overline{q_p}) = \ker(\widehat{q_p})\oplus \ker(\widehat{q_p})^\perp
\xrightarrow{\partial_p=
\left(\begin{smallmatrix} \overline{\partial_p} &\gamma'\\
 0 & pr^{-1}\delta_p\widehat{q_p}\end{smallmatrix}\right)}
 \overline{\im(c^p)}\oplus (\im(c^p)^\perp\cap \ker(c_{p+1}))=
 \ker(c_{p+1}).
\end{multline*}
Because of $\norm{\gamma'}\le \norm{\partial_p}$
Lemma \ref{matrix_spec_dens}.3 implies
\[ \overline{F}(\partial_p,\lambda) \le
\overline{F}(\overline{\partial_p},\lambda^{1/2}) + \overline{F}(pr^{-1}
 \delta_p\widehat{q_p},(4 + 2\norm{\partial_p})
 \lambda^{1/2})\quad\text{for }
 \lambda<(4 + 2\norm{\partial_p})^{-2}.\]
Note that $\widehat{q_p}^{-1}= q_p^{-1}\circ (pr:\ker e^p\to H^p(E))^{-1}$,
therefore $\norm{\widehat{q_p}^{-1}}\le \norm{q_p^{-1}}$. Moreover,
$\partial_p= j_{p+1}^{-1}\circ\overline{d^p}$, hence $\norm{\partial_p}\le
\norm{j_{p+1}^{-1}}\cdot\norm{\overline{d^p}}=
\norm{j_{p+1}^{-1}}\cdot\norm{d^p}$. Since a non-trivial
projection and its inverse always have norm $1$,
we conclude
from Lemma \ref{basic_barF} 4.) and 5.) that
for all $\lambda < (4 + 2\norm{j_{p+1}^{-1}}\norm{d^p})^{-2}$
\begin{eqnarray}
\overline{F}(\partial_p,\lambda) & \le &
 \overline{F}(\overline{\partial_p},\lambda^{1/2}) +
\overline{F}(\delta_p,\norm{q_p^{-1}}
(4 + 2\norm{j_{p+1}^{-1}}\norm{d^p})\lambda^{1/2}). \label{eqn 4.12}
\end{eqnarray}
It remains to identify $\overline{\partial_p}$.
Lott-L\"uck \cite[page 29]{Lott-Lueck (1995)} define
maps so that the following diagram is commutative:
\[\begin{CD}
 C^p @= C^p @>{\widetilde{\overline{j_p}}}>> \ker{\widehat{q_p}}\\
@VV{j_p}V @VV{\overline{j_p}}V @VV{incl}V\\
D^p @>>>D^p/\overline{\im(d_{p-1})}
@>{\id}>> D^p/\overline{\im(d_{p-1})}
\end{CD}\]
and show that $ c^p =
\overline{\partial_p}\circ \widetilde{\overline{j_p}}$
and that $\widetilde{\overline{j_p}}$ has dense image and is left
Fredholm. Because of $\norm{\widetilde{\overline{j_p}}}
= \norm{\overline{j_p}}\le\norm{j_p}$ Lemma
\ref{basic_barF}.2 implies
\begin{eqnarray}
\overline{F}(\overline{\partial_p},\lambda) & \le &
\overline{F}_p(C^*,\norm{j_p}\lambda). \label{eqn 4.13}
\end{eqnarray}
Now Theorem \ref{short_exact} follows from
\eqref{eqn 4.10}, \eqref{eqn 4.11}, \eqref{eqn 4.12} and \eqref{eqn 4.13}.
\qed



\typeout{--------------------- section 5 -----------------------}

\section{Sobolev- and $L^2$-complexes}
\label{Sobolev- and L^2-complexes}

In this section we show how
the study of the spectral density function of the
Laplacian with absolute boundary conditions,
considered as an unbounded operator on $L^2$, can be translated to
the study of the spectral density functions of Sobolev de Rham
complexes without any boundary conditions.

As intermediate steps we study an
$L^2$-de Rham complex with absolute
boundary conditions, then a Sobolev
complex with absolute boundary conditions,
and in a last step the Sobolev complex
without boundary conditions. We need
the last one, because only here,
an exact Mayer-Vietoris sequence can be
obtained. Efremov \cite{Efremov (1991)}
and Lott-L\"uck \cite{Lott-Lueck (1995)}
use the same reductions. We repeat and
refine their arguments, because we need
more precise information on the spectral
density functions than they do.

\begin{definition} \label{Sobolev spaces}
Let $N$ be a complete $m$-dimensional
Riemannian manifold. If $N$ has no boundary,
define for each integer $p \ge 0$ and each real number $s \ge 0$
the {\em Sobolev norm} $\abs{\;}_{H^s}$
on the space $C^{\infty}_0(\Lambda^p(N))$
of smooth $p$-forms with compact support by
\[\abs{\omega}_{H^s} ~ := ~
\abs{(1 + \Delta_p)^{s/2}\omega}_{L^2}.\]
If $N$ has boundary, define (compare \cite[page 363]{Taylor (1996a)})
for integers $p,s \ge 0$ the {\em Sobolev norm}
$\abs{\;}_{H^s}$ on the space $C^{\infty}_0(\Lambda^p(N))$ inductively
by
\begin{eqnarray*}
\abs{\omega}_{H^0} & = & \abs{\omega}_{L^2};
\\
\abs{\omega}_{H^{s+1}}^2 & = & \abs{\omega}^2_{H^s}+
\abs{d^p\omega}_{H^{s}}^2 + \abs{\delta^p\omega}_{H^{s}}^2 +
\abs{i^*(*\omega)}_{H^{s+1/2}}^2,
\end{eqnarray*}
where $i: \partial N \longrightarrow N$ is the inclusion of the boundary
and $\delta^p = (-1)^{p(m-p)}) * d^{m-p} *$ is the adjoint of $d^p$
with respect to the $L^2$-inner product. The {\em Sobolev space}
$H^s(\Lambda^p(N))$ is the completion of $C^{\infty}_0(\Lambda^p(N))$
with respect to $\abs{\,}_{H^s}$. \qed
\end{definition}

Next we introduce the various relevant Sobolev and $L^2$-chain complexes.

\begin{definition} \label{Sobolev cochain complexes}
Let $N$ be a complete $m$-dimensional
Riemannian manifold. Define its {\em Sobolev cochain complex} $D_p^*[N]$
which is concentrated in dimensions $p-1$, $p$ and $p+1$ by
\[\ldots \to 0 \to H^2(\Lambda^{p-1}(N))\xrightarrow{d^{p-1}}
H^1(\Lambda^p(N))\xrightarrow{d^p}
L^2(\Lambda^{p+1}(N))\to 0 \to \ldots .\]

Define the {\em Sobolev cochain complex
with absolute boundary conditions}  $D_{p,abs}^*[N]$
which is again concentrated in dimensions $p-1$, $p$, $p+1$ by
\[ \ldots \to 0 \to H^2_{abs}(\Lambda^{p-1}(N))\xrightarrow{d}
H^1_{abs}(\Lambda^p(N))\xrightarrow{d}
L^2(\Lambda^{p+1}(N))\to 0 \to \ldots, \]
where
\begin{eqnarray*}
 H^1_{abs}(\Lambda^{p-1}(N))& := &
\{\omega\in H^1(\Lambda^{p-1}(N))\mid i^*(*\omega)=0\}
~ \subset H^1(\Lambda^{p-1}(N));\\
H^2_{abs}(\Lambda^{p-1}(N))& := &
\{\omega \in H^2(\Lambda^{p-1}(N))\mid
i^*(*\omega)=0=i^*(*d\omega)\}~ \subset H^2(\Lambda^{p-1}(N)).
\end{eqnarray*}
The inclusions induce a cochain map and we will use
on $D_{p,abs}^*[N]$ the norm induced from  $D_p^*[N]$.

Define the {\em $L^2$-cochain complex} $L^*_p[N]$
which is again concentrated in dimensions $p-1$, $p$, $p+1$ by
\[\ldots \to 0 \to L^2(\Lambda^{p-1}(N))\xrightarrow{d^{p-1}}
L^2(\Lambda^p(N))\xrightarrow{d^p}
L^2(\Lambda^{p+1}(N))\to 0 \to \ldots,\]
where $d^p$ is the closure of
the operator with domain $C^\infty_0(\Lambda^pN)$. \qed
\end{definition}

Notice that the differentials in
$D_{p}^*[N]$ and $D_{p,abs}^*[N]$ are bounded operators,
what is not true for $L^*_p[N]$.

\subsection{Sobolev complexes with and
without boundary conditions}
\subsubsection{Arbitrary manifolds}

Let $N$ be a compact Riemannian manifold possibly with boundary.
The inclusion induces a cochain map
$i^*: D^*_{p,abs}[\widetilde{N}]\to D^*_p[\widetilde{N}]$
(see Definition \ref{Sobolev cochain complexes}).
Given a geodesic collar of  width $w>0$, Lott-L\"uck
\cite[Definition 5.4 and Lemma 5.5]{Lott-Lueck (1995)}
define maps
\[ K_s^p: H^s(\Lambda^p(\widetilde{N}))\to
H^{s+1}(\Lambda^{p-1}(\widetilde{N}))
\qquad s =0,1,2, \]
which induce maps
\[ K_{s,abs}^p: H^s_{abs}(\Lambda^p(\widetilde{N}))\to
H^{s+1}_{abs}(\Lambda^{p-1}(\widetilde{N}))
\qquad s =0,1\]
and show that they have the following properties
\begin{lemma}\label{K_prop}
\begin{enumerate}
\item The maps $K_s^p$ are bounded;
\item $K_s^p\omega$ depends only on the restriction of $\omega$
 to the interior of the collar of width $w$, and $\supp K_s^p\omega$ is
contained in this collar;
\item The maps
\[
\begin{split}
  j_p^{p-1}&:= 1-d^{p-2}K_2^{p-1}- K_1^{p-1}d^{p-1};\\
  j_p^p &:= 1- d^{p-1}K_1^p - K_0^{p+1}d^p;\\
  j_p^{p+1}& := 1- d^p K_0^{p+1};
\end{split} \]
constitute a chain map
\[ j_p^*: D^*_p[\widetilde{N}]\to D^*_{p,abs}[\widetilde{N}];\]
\item The inclusion $i^*$ and the map $j^*$ induce (inverse)
homotopy equivalences
$\overline{i}_p^*: \overline{D}_p^* \to \overline{D}_{abs,p}^* $ and
$\overline{j}_p^*: \overline{D}_{abs,p}^* \to \overline{D}_p^*$
between the reduced complexes
\begin{eqnarray*}
 \overline{D}^*_p & := &
\ldots \to 0 \to D^p_p[\widetilde{N}]/\clos(\im d^{p-1})
\to D^{p+1}_p[\widetilde{N}]  \to 0 \to \ldots \\
\widetilde{D}^*_{p,abs} &:= &
\ldots \to 0 \to D^p_{abs,p}[\widetilde{N}]/\clos(\im d^{p-1})
\to D^{p+1}_{p,abs}[\widetilde{N}]  \to 0 \to \ldots
\end{eqnarray*}
The corresponding homotopies to
the identity are induced from
$K^*_*\circ i^*$ and $i^*\circ K^*_*$. \qed
\end{enumerate}
\end{lemma}

\subsubsection{Application to hyperbolic manifolds}

We use these results to compare the
spectral density function of
$D^*_p[\widetilde{M_R}]$ and $D^*_{p,abs}[\widetilde{M_R}]$.
Note that for $R\ge 1$, $\widetilde{M_R}$
has a geodesic collar of width $1/3$,
and all these collars are isometric to the one of $\widetilde{M_1}$.
In particular we get
\begin{eqnarray}
\norm{K^p_s[\widetilde{M_R}]}& \le & C \hspace{10mm} s = 0,1,2,~ p \ge 0;
\label{op_est1}
\\
\norm{j^q_p[\widetilde{M_R}]} & \le & C
\hspace{10mm} p \ge 0, ~ q = p-1,p,p+1;
\label{op_est2}
\end{eqnarray}
with a constant $C$ not depending on $R$
since the maps $K^p_s$ involve only the collar of the
boundary. Now we can use the following theorem of Gromov-Shubin
\cite[Proposition 4.1]{Gromov-Shubin (1991)}.
\begin{theorem} \label{Gromov-Shubin}
  Let $C^*$ and $D^*$ be cochain complexes of Hilbert
$\Neumann{N}(\Gamma)$-modules with not necessarily bounded differential,
 $f^*:C^*\to D^*$
and $g^*:D^*\to C^*$ bounded
homotopy equivalences and $T^*:C^*\to C^{*-1}$
a homotopy between $g^*f^*$ and
$\id$. Then for the spectral density functions the following holds
\[ F_p(C^*,\lambda) \le F_p(D^*, \norm{f_{p+1}}^2\norm{g_p}^2\lambda)
 \qquad\forall \lambda< (2\norm{T_{p+1}})^{-2}. \qed \]
\end{theorem}

\begin{proposition}\label{Fohne}
We find constants $C_1,C_2>0$, independent of $R$, so that
$$\begin{array}{lclcll}
F_p(D^*_{p,abs}[\widetilde{M_R}],C_1^{-1}\lambda)
& \le &
F_p(D^*_p[\widetilde{M_R}],\lambda)
& \le & F_p(D^*_{p,abs}[\widetilde{M_R}], C_1\lambda)
& \qquad\forall
\lambda\le C_2;
\\
F_p(D^*_{p,abs}[\widetilde{T_R}],C_1^{-1}\lambda)
& \le &
F_p(D^*_p[\widetilde{T_R}],\lambda)
& \le & F_p(D^*_{p,abs}[\widetilde{T_R}], C_1\lambda)
& \qquad\forall
\lambda\le C_2.
\end{array}$$
\end{proposition}
\proof
For $\widetilde{M_R}$ this is a direct consequence of
Theorem \ref{Gromov-Shubin}
applied to the homotopy equivalence in Lemma \ref{K_prop}.4,
the estimates \eqref{op_est1} and \eqref{op_est2} and of the
fact, that the $p$-spectral density
function of $D^*_p[\widetilde{N}]$ and $D^*_{p,abs}[\widetilde{N}]$,
respectively, coincide by definition with the on
of $\overline{D}^*_p[\widetilde{N}]$ and
$\overline{D}^*_{p,abs}[\widetilde{N}]$,
respectively. The case of $\widetilde{T_R}$
is completely analogous.\qed

\subsection{$L^2$-complexes and
Sobolev complexes with boundary conditions}

To compare spectral density functions of
Sobolev complexes and $L^2$-de Rham
complexes, we need the following formula
for these functions:
\begin{lemma}\label{compute}
Suppose $(E^*,d^*)$ is a Hilbert-$\Neumann{A}$
cochain complex. Here,
we use the broader definition of
Gromov-Shubin \cite[section 4]{Gromov-Shubin (1991)},
where $d^*$ are closed, but not necessarily
bounded operators. Then
\[ F_p(E^*,\lambda) = \sup_{L\in S^p_\lambda} \dim_\Gamma L, \]
where $S^p_\lambda$ is the set of all closed
$\Gamma$-invariant subspaces $L$ of
$\ker(d^p)^\perp\cap\mathcal{D}(d^p)\subset E^p$ so that
\[ \abs{d^p x}\le \lambda\cdot \abs{x}\qquad\forall x\in L. \]
\end{lemma}
\begin{proof}
The proof in \cite[Lemma 1.5]{Lott-Lueck (1995)},
where the proposition is stated
only for bounded $d^*$, works without modifications also for
unbounded operators.
\end{proof}

To apply this theorem it is necessary
to compute the orthogonal complement of
the kernel of the differential for the considered complexes.

\begin{lemma} \label{orthogonal complement of kernel}
Let $N$ be a compact Riemannian manifold possibly with boundary. Then
\begin{eqnarray*}
\ker\left(d^p: H^1_{abs}(\Lambda^p(\widetilde{N}))
\to L^2(\Lambda^{p+1}(\widetilde{N}))\right)^{\perp_{H^1}}
& = &
H^1_{abs}(\Lambda^p(\widetilde{N})) \cap
\overline{\delta C^\infty_{abs}(\Lambda^{p+1}(\widetilde{N}))}^{L^2},
 \end{eqnarray*}
where $C^\infty_{abs}(\Lambda^{p+1}(\widetilde{N})):=
\{\omega\in C^\infty_0(\Lambda^{p+1}(\widetilde{N}))\mid  i^*(*\omega)=0\}$
with $i: \partial \widetilde{N} \hookrightarrow \widetilde{N}$
the inclusion.
\end{lemma}
\proof  $H^1_{abs}(\Lambda^p(\widetilde{N})) \cap
\overline{\delta C^\infty_{abs}(\Lambda^{p+1}(\widetilde{N}))}^{L^2}
 ~ \subset ~
\ker\left(d^p: H^1_{abs}(\Lambda^p(\widetilde{N}))
\to L^2(\Lambda^{p+1}(\widetilde{N}))\right)^\perp$:

 Let $x\in H^1_{abs}\cap
\overline{\delta C^\infty_{abs}}^{L^2}$ with $x=\lim_{L^2}\delta x_n$,
 $y\in H^1_{abs}$ with $dy=0$. Then
\[\begin{split}
  (x,y)_{H^1} &=(x,y)_{L^2}+(dx,\underbrace{dy}_{=0})_{L^2}+(\delta x,\delta y)_{L^2}
 + (i^*(*x),\underbrace{i^*(*y)}_{=0})_{H^{1/2}}\\
(x,y)_{L^2} &= \lim(\delta x_n,y)_{L^2}=\lim (x_n,dy)_{L^2}\pm \int_{\boundary
\tilde N}y\wedge *x_n =0
\end{split}\]
If $z\in C^\infty_0(\tilde N-\boundary\tilde N)$ then
\[ (\delta x,z)=(x,dz)+\int_{\boundary\tilde N} z\wedge *x =\lim(\delta x_n,dy)
=\lim(\underbrace{\delta\delta}_{=0}x_n,z) =0 \]
The set of these $z$ is dense in $L^2$ hence $\delta x=0$.

It remains to show the opposite inclusion.
Put
$$\mathcal{H}^p ~ = ~
\{x\in H^\infty(\Lambda^p(\widetilde{N}))\mid
dx=0=\delta x\text{ and }i^*(*x)=0 \}.$$
(Because of elliptic regularity, we could replace $H^{\infty}$ by
$C^{\infty} \cap L^2$).
There is the orthogonal decomposition \cite[Theorem 5.10]{Schick (1996)} or
\cite{Burghelea-Friedlander-Kappeler (1996a)}:
\begin{equation}\label{L2split}
L^2(\Lambda^p(\widetilde{N})) = \mathcal{H}^p\oplus
\overline{d C^\infty_0(\Lambda^{p-1}(\widetilde{N}))}\oplus
\overline{\delta C^\infty_{abs}(\Lambda^{p+1}(\widetilde{N}))}.
\end{equation}
Hence it suffices to show that
$y \in \ker\left(d^p: H^1_{abs}(\Lambda^p(\widetilde{N}))
\to L^2(\Lambda^{p+1}(\widetilde{N}))\right)^{\perp_{H^1}}$
 is orthogonal to both $\mathcal{H}^p$
and $\overline{d C^\infty_0(\Lambda^{p-1}(\widetilde{N}))}$
in $L^2(\Lambda^p(\widetilde{N}))$.
The first claim follows from
\[ 0 ~ = ~ (x,y)_{H^1} ~ = ~ (x,y)_{L^2}
\quad\forall x\in \mathcal{H}^p, \]
and for the second claim it suffices
to prove $\delta y=0$. We have
\[ (dx,y)_{H^1}=0\quad\forall x\in H^\infty(\Lambda^{p-1}(\widetilde{N}))
\text{ with $i^*(*dx)=0$}. \]
Since $i^*(*y)=0$ this implies
\[
\begin{split}
   0 &=(dx,y)_{L^2}+(\delta dx,\delta y)_{L^2} +
 (\underbrace{dd\delta x}_{=0},y)_{L^2}\\\
     &=(x,\delta y)_{L^2} +
       (\delta dx,\delta y)_{L^2}+ (d\delta x,\delta y)_{L^2}\\
     &= ((1+\Delta_{p-1}) x,\delta y)_{L^2}.
\end{split} \]
Now
$$(1+\Delta_{p-1}):H^\infty_{abs}(\Lambda^{p-1}(\widetilde{N}))=
\{x\in H^\infty(\Lambda^{p-1}(\widetilde{N}));\;i^*(*x)=0=i^*(*dx)\}
\to L^2(\Lambda^{p-1}(\widetilde{N}))$$
has dense image (compare \cite[Theorem 5.19]{Schick (1996)}
and  therefore $\delta y=0$. This finishes the proof
of Lemma \ref{orthogonal complement of kernel}. \qed

\begin{lemma}\label{ort=}
Let $N$ be a compact Riemannian manifold possibly with boundary
and $L^*_p[\widetilde{N}]$ be the cochain complex introduced
in Definition \ref{Sobolev cochain complexes}. Then
\[ \mathcal{D}(d)\cap\ker\left(d^p: L^2(\Lambda^p(\widetilde{N}))\to
L^2(\Lambda^{p+1}(\widetilde{N}))\right)^{\perp_{L^2}}
~ = ~
H^1_{abs}(\Lambda^p(\widetilde{N}))\cap
\overline{\delta C^\infty_{abs}(\Lambda^{p+1}(\widetilde{N}))}^{L^2} \]
with $C^\infty_{abs}$
as in Lemma \ref{orthogonal complement of kernel}.
\end{lemma}

\begin{proof}
First remember that $d+\delta:C^\infty(\Lambda^\bullet\widetilde{N})\to
C^\infty(\Lambda^\bullet(\widetilde{N}))$ with either absolute
(i.e.\ $i^*(*\omega)=0$) or relative (i.e.\ $i^*(\omega)=0$)
 boundary conditions are formally self adjoint
elliptic boundary value problems. We will establish that
$H^1_{abs}(\Lambda^p(\widetilde{N}))\cap
\overline{\delta C^\infty_{abs}(\Lambda^{p+1}(\widetilde{N}))}^{L^2}$ is
perpendicular to $\ker(d^p)$, and that any form which
is perpendicular to
$\ker(d^p)$ and is contained in the domain of $d$ lies in
$H^1_{abs}(\Lambda^p(\widetilde{N}))\cap
\overline{\delta C^\infty_{abs}(\Lambda^{p+1}(\widetilde{N}))}$.

For the first
statement take $x \in C^\infty_{abs}(\Lambda^{p+1}(\widetilde{N}))$
and $y\in\ker d^p$. In particular, a
sequence $y_n\in C_0^\infty(\Lambda^{p+1}(\widetilde{N}))$ exists with
$y_n\xrightarrow{L^2}y$ and $d y_n
\xrightarrow{L^2} 0$. Then
\[ (y,\delta x)_{L^2} =\lim (y_n,\delta x)_{L^2}\stackrel{i^*(*x)=0}{=}
\lim(dy_n,x)_{L^2} =0. \]
Therefore, $H^1_{abs}(\Lambda^p(\widetilde{N}))\cap
\overline{\delta C^\infty_{abs}(\Lambda^{p+1}(\widetilde{N}))}^{L^2}$
and $\ker(d^p)$ are perpendicular. For
the second statement, suppose $x\in\mathcal{D}(d)$
is perpendicular to $\ker d^p$.
Choose $x_n\in C^\infty_0(\Lambda^p(\widetilde{N}))$ with
$x_n\xrightarrow{L^2}x$
and $dx_n\to dx$. For
every $y\in C^\infty_{abs}(\Lambda^{\bullet}(\widetilde{N}))$ we have
\[
\begin{split}
  ((d+\delta)y,x)_{L^2} \stackrel{dy\in\ker d}{=}&(\delta y,x)_{L^2}=
  \lim_{n \to \infty}(\delta y,x_n)_{L^2}\\
 \stackrel{i^*(*y)=0}{=} & \lim_{n \to \infty}(y, dx_n)_{L^2} =
(y,dx)_{L^2}
\end{split} \]
Adjoint elliptic regularity
\cite[Lemma 4.19, Corollary 4.22]{Schick (1996)}
shows that $x\in H^1_{loc}(\Lambda^p(\widetilde{N}))$. Then
 $(\delta x,y)_{L^2}
=(x,dy)=0$ holds since $\forall y\in C^\infty_0(\tilde N-\boundary\tilde N)$
 $dy\in\ker d$. We conclude $ \delta x=0$. This means that
$\forall y\in C^\infty_0(\tilde N)$
\[ 0\stackrel{dy\in\ker d}{=} (dy,x)=(y,\underbrace{\delta x}_{=0})
\pm\int_{\boundary\tilde N} y\wedge*x \]
Since $\{i^*(y)\}$ is dense in $L^2(\boundary \tilde N)$,
 $i^*(*x)=0$, i.e.\
$x\in H^1_{abs}(\Lambda^p(\widetilde{N}))$. The $L^2$-splitting
\eqref{L2split} implies $x\in\overline{\delta C^\infty_{abs}}$.
\end{proof}

Now we can compare the spectral density
functions of $L^*_p$  and
$D^*_{abs,p}$.
\begin{proposition}\label{Fabs}
Let $N$ be any compact manifold with
boundary, $\Gamma:=\pi_1(N)$. Then
\[
\begin{split}
  F_p(L^*_p[\widetilde{N}],\lambda)
 &\le F_p(D^*_{abs,p}[\widetilde{N}], \lambda)\qquad\forall \lambda;
 \\
  F_p(D^*_{abs,p}[\widetilde{N}],\lambda)
&\le F_p(L^*_p[\widetilde{N}],\sqrt{2}\lambda)\qquad\forall\lambda\le
    1/\sqrt{2}.
\end{split}\]
\end{proposition}
\proof
We will use Lemma \ref{compute}.

We start with the first inequality.
Let $L\subset \ker\left(d^p :L^2(\Lambda^p(\widetilde{N})) \to
L^2(\Lambda^{p+1}(\widetilde{N}))\right)^{\perp_{L^2}}
\cap \mathcal{D}(d^p)$
be a closed $\Gamma$-invariant subspace with
$\abs{dx}_{L^2}\le\lambda\abs{x}_{L^2}$
for all $x\in L$. Hence we get
$\abs{dx}_{L^2}\le\lambda\abs{x}_{H^1}$
for all $x\in L$.
Moreover, $L\subset\ker\left(d^p :H^1(\Lambda^p(\widetilde{N}))
\to L^2(\Lambda^{p+1}(\widetilde{N}))\right)^{\perp_{H^1}}$ since
the two orthogonal complements
are equal by Lemma \ref{ort=} and Lemma \ref{orthogonal complement of kernel}.
$L$ is closed also with respect to the $H^1$-topology because this is finer
than the $L^2$-topology.

By Lemma \ref{compute}, $F_p(L^*_p(\widetilde{N}),\lambda)=
\sup_L \dim_\Gamma L$, where the supremum
is over all such $L$. We have just seen
that for the computation of
$F_p(D^*_{abs}(\widetilde{N}),\lambda)$ we have to
take the supremum over a larger set. This
implies the first inequality. It remains to prove the second.

Let $\lambda\le 1/\sqrt{2}$ and let $L\subset
H^1_{abs}(\Lambda^p(\widetilde{N}))$ be
a closed $\Gamma$-invariant subspace with
$L\perp_{H^1} \ker(d^p :H^1_{abs}(\Lambda^p(\widetilde{N}))\to
L^2(\Lambda^{p+1}(\widetilde{N})))$ and
$\abs{dx}_{L^2}\le\lambda\abs{x}_{H^1}$ for all
$x\in L$. Then Lemma \ref{ort=} implies for all $x \in L$ that
\begin{eqnarray}
\abs{x}_{H^1}^2  & = & \abs{x}_{L^2}^2 + \abs{dx}_{L^2}^2.
\label{eqn 6.1}
\end{eqnarray}
\begin{equation}
\implies\quad \abs{x}_{L^2} ~ \le ~ \abs{x}_{H^1} ~ \le ~
(1 - \lambda^2)^{-1/2}\cdot  \abs{x}_{L^2}.
\label{eqn 6.2}
\end{equation}
Equation \eqref{eqn 6.2} says that on $L$ the $L^2$-norm and the
$H^1$-norm are equivalent, so that $L$ is a closed subspace
of $L^2(\Lambda^p(\widetilde{N}))$. We conclude from \eqref{eqn 6.1}
\[\begin{split}
&\abs{dx}_{L^2}^2
\le \lambda^2 \abs{x}_{H^1}^2 \le
\lambda^2\abs{x}^2_{L^2}+\lambda^2\abs{dx}_{L^2}^2
\le
\lambda^2 \abs{x}^2_{L^2} + \frac{1}{2}\abs{dx}^2_{L^2}\\
\implies& \abs{dx}_{L^2}\le \sqrt{2}\lambda\abs{x}_{L^2}\\
\end{split}\]
Now the second inequality follows as above. \qed

\subsection{$L^2$-complexes and the Laplacian}

Let $N$ be a compact manifold
with boundary as above. In the last paragraph,
we studied the unbounded operator $d$ on $\ker(d)^\perp\subset
L^2(\widetilde{N})$ with domain
$H^1_{abs}\cap \overline{\delta C^\infty_{abs}}$.
Obviously, this coincides with
the unbounded operator $d+\delta$ with domain
$H^1_{abs}$, restricted to
$\overline{\delta C^\infty_{abs}}$. This boundary
value problem is elliptic,
 hence $(d+\delta)$ with domain
$H^1_{abs}$ is self adjoint in the Hilbert space sense
(\cite[Theorem 4.25]{Schick (1996)}).
The adjoint of $d$ restricted to
$ \overline{\delta C^\infty_{abs}}$
is therefore $d+\delta$ composed with
projection onto $ \overline{\delta C^\infty_{abs}}$.
The square of $d+\delta$ is just
\[\Delta=(d+\delta)^2\qquad\text{with domain }
H^2_{abs}=\{\omega\in H^2;\;
 i^*(*\omega)=0=i^*(*d\omega)\} \]
This is exactly the operator we have to study.
We want to compare its spectral density function
with the one  of $L_p^*(\widetilde{N})$.
Recall that $\Delta_p^{\perp}$ is the operator from
the orthogonal complement of the kernel
of $\Delta_p$ to itself obtained by restriction. Note that
the splitting \eqref{L2split} of $L^2$
induces a splitting of the Laplacian:
\begin{equation}
  \label{splitting of Laplacian in summands}
 \Delta_p^{\perp}[\widetilde{N}] ~ = ~
\delta^{p+1} d^p|_{\overline{\delta
C^\infty_{abs}(\Lambda^{p+1}(\widetilde{N}))}}\oplus
d^{p-1}\delta^p|_{\overline{dC^\infty_0(\Lambda^{p-1}(\widetilde{N}))}}. 
\end{equation}

\begin{lemma} \label{computations of dd^* and d^*d}
We have
\begin{eqnarray*}
(\delta^{p+1} d^p)|_{\overline{\delta
C^\infty_{abs}(\Lambda^{p+1}(\widetilde{N}))}}
& = &
(d^p|_{\overline{\delta C^\infty_{abs}(\Lambda^{p+1}(\widetilde{N}))}})^*
(d^p|_{\overline{\delta C^\infty_{abs}(\Lambda^{p+1}(\widetilde{N}))}});
\\
(d^{p-1}\delta^p)|_{\overline{d C^\infty_0(\Lambda^{p-1}(\widetilde{N}))}}
& = &
(d^{p-1}|_{\overline{\delta
C^\infty_{abs}(\Lambda^{p-1}(\widetilde{N}))}})
(d^{p-1}|_{\overline{\delta
C^\infty_{abs}(\Lambda^{p-1}(\widetilde{N}))}})^*.
\end{eqnarray*}
Here
$d^p|_{\overline{\delta C^\infty_{abs}(\Lambda^{p+1}(\widetilde{N}))}}$ is
the unbounded operator on the subspace
$\overline{\delta C^\infty_{abs}(\Lambda^{p+1}(\widetilde{N}))}$ of
$L^2(\Lambda^p(\widetilde{N}))$ with range
$\overline{d C^\infty_0(\Lambda^p(\widetilde{N}))}$ and with domain
$H^1_{abs}(\Lambda^p(\widetilde{N}))\cap
\overline{\delta C^\infty_{abs}(\Lambda^{p+1}(\widetilde{N}))}$.
\end{lemma}
\proof
 We first prove that the Hilbert space adjoint $d^*$
of $d|_{\overline{\delta C^\infty_{abs}}}: \overline{
\delta C^\infty_{abs}}\to
\overline{d C^\infty_0}$ is $\delta$ with domain
$\overline{d C^\infty_0}\cap H^1_{abs}$.

If $x\in\mathcal{D}(d^*)\subset
\overline{ d C^\infty_0}$ $\implies$
$x\in\mathcal{D}(d)$ and $dx=0$.

 Moreover,
for arbitrary $y\in H^1_{abs}\cap\overline{\delta C^\infty_{abs}}$ we have
$\delta y=0$. Therefore
\[ ((d+\delta)y,x)_{L^2}=(dy,x)_{L^2} = (y,d^*x)_{L^2}. \]
If $y\in H^1_{abs}\cap\overline{d C^\infty_0}$ then
\[ ((d+\delta)y,x))_{L^2} = (\delta y,x))_{L^2}=0 \]
because $\delta H^1_{abs}\perp d C^\infty_0$.
But also $(y,d^*x))_{L^2}=0$ because
$d^*x\in \overline{\delta C^\infty_{abs}}$
and $d C^\infty_0\perp
\delta C^\infty_{abs}$. Adjoint elliptic
regularity \cite[Lemma 4.19]{Schick (1996)}
 implies
now that $x\in H^1_{loc}$. We have to show $x\in H^1_{abs}$,
 i.e.~$dx,\delta x\in L^2$
and $i^*(*x)=0$. Now
\[ dx=0\in L^2;\qquad \delta x\stackrel{dx=0}{=}(d+\delta)x = d^*x \in L^2. \]
\[ \left(
(x,d\delta y)=(d^*x,\delta y)=(\delta x,\delta y)= (x,d\delta y)\pm\int_{
\boundary\tilde N}\delta y\wedge*x\quad\forall y\in C^\infty_0(\tilde N)
\right)
\quad\implies i^*(*x)=0 \]

Clearly $\delta d|_{\mathcal{D}(d^*d)}=\Delta=\Delta^\perp$. To prove the
lemma
it remains to show
that the domains coincide, i.e.\ that
\[ \mathcal{D}(\Delta)\cap \overline{\delta C^\infty_{abs}} =
  \mathcal{D}(d^* (d|_{\overline{\delta C^\infty_{abs}}})) \]
Now $\mathcal{D}(\Delta)= H^2_{abs}=\{x\in H^2;\;i^*(*x)=0=i^*(*dx)\}
\subset\mathcal{D}(d^*d)$. 

If, on the other hand, $x\in\mathcal{D}(d^*d|)$ then
$x\in H^1_{abs}\cap\overline{
\delta C^\infty_{abs}}$ and $dx\in H^1_{abs}$. I.e.\
$(d+\delta)x\in H^1$ because  $\delta x=0$. Since also $i^*(*x)=0$,
by elliptic regularity
$x\in H^2$. The boundary conditions are
fulfilled therefore $x\in H^2_{abs}=
\mathcal{D}(\Delta)$.

The proof for $(d\delta)|_{\overline{d C^\infty_0}}$ is similar.
\qed

\begin{proposition} \label{comparision of Laplace and L}
Let $N$ be a compact Riemannian manifold possibly with boundary.
Then
\[F(\Delta_p^{\perp}[\widetilde{N}],\sqrt{\lambda}) ~ = ~
F_p(L^*_p[\widetilde{N}],\lambda) + F_{p-1}(L^*_p[\widetilde{N}],\lambda).
\]
\end{proposition}
\proof This follows from 
 \eqref{splitting of Laplacian in summands}, Lemma \ref{matrix_spec_dens} 1.),
Lemma \ref{computations of dd^* and d^*d}, Lemma \ref{basic_F} 6.) and
Lemma \ref{F and adjoint}
by the following calculation.
\begin{eqnarray*}
F(\Delta_p^{\perp}[\widetilde{N}],\sqrt{\lambda})
& = & F(\delta^{p+1} d^p|_{\overline{\delta
C^\infty_{abs}(\Lambda^{p+1}(\widetilde{N}))}}, \sqrt{\lambda}) +
F(d^{p-1}\delta^p|_{\overline{dC^\infty_0(\Lambda^{p-1}(\widetilde{N}))}},
\sqrt{\lambda})
\\
& = &
F(d^p|_{\overline{\delta C^\infty_{abs}(\Lambda^{p+1}(\widetilde{N}))}}),
\lambda) +
F((d^{p-1}|_{\overline{\delta
C^\infty_{abs}(\Lambda^{p-1}(\widetilde{N}))}},\lambda)
\\
& = &
F_p(L^*_p[\widetilde{N}],\lambda) + F_{p-1}(L^*_{p-1}[\widetilde{N}],\lambda).
\qed \end{eqnarray*}



\typeout{--------------------- section 6 -----------------------}

\section{Spectral density functions for $M_R$}
\label{Spectral density functions for M_R}

In this section, we estimate
the spectral density function of the
Laplacian on $\widetilde{M_R}$ independently of $R$
and finish the proof of the main Theorem \ref{Tor_converg}.

The following sequence of Hilbert cochain complexes is exact
\cite[Lemma 5.14]{Lott-Lueck (1995)}
\begin{equation}\label{ses}
 0\to D_p^*[\widetilde{M}] \xrightarrow{j}
 D_p^*[\widetilde{M_R}] \oplus  D_p^*[\widetilde{E_{R-1}}]
\xrightarrow{q} D_p^*[\widetilde{T_{R-1}}] \to 0
\end{equation}
with $j(\omega):=i^*_{M_R}\omega\oplus i^*_{E_{R-1}}\omega$ and
$q(\omega_1 \oplus\omega_2):=i^*_{T_{R-1}}\omega_1 - i^*_{T_{R-1}}\omega_2$,
where we use Notation \ref{notation for the ends} and
$i_X$ are inclusion maps. We are interested in the spectral
density function of $D_p^*[\widetilde{M_R}]$ at $* = p$. We get from
Lemma \ref{F and block matrices}.1
\begin{eqnarray}
F_p(D^*_p[\widetilde{M_R}],\lambda)  & \le &
F_p(D^*_p[\widetilde{M_R}] \oplus D^*_p[\widetilde{E_{R-1}}],\lambda).
\label{eqn 5.0}
\end{eqnarray}
The latter can be estimated in
terms of the spectral densities of $D_p^*[\widetilde{M}]$
and $D_p^*[\widetilde{T_{R-1}}]$ by Theorem \ref{short_exact}.
To apply this theorem, we have to check the Fredholm condition,
compute the connecting homomorphism
$\delta^p$ of the long weakly exact $L^2$-cohomology sequence
of our short exact sequence and
the constants appearing in the statement of
the Theorem \ref{short_exact}.

\begin{lemma} \label{delta is trivial}
\begin{enumerate}
\item We get for all $R \ge 1$ and $\lambda \ge 0$
\begin{eqnarray}
F(\Delta^{\perp}_p[\widetilde{T_{R-1}}],\lambda)
& \le
F(\Delta^{\perp}_p[\widetilde{T_0}],\lambda);
\end{eqnarray}

\item There are positive numbers
$\alpha > 0$  and $\epsilon > 0$ such that for $\lambda \le \epsilon$
\begin{eqnarray*}
F(\Delta_p^{\perp}[\widetilde{M}],\lambda) & \le & \lambda^{\alpha};
\\
F(\Delta_p^{\perp}[\widetilde{T_0}],\lambda) & \le & \lambda^{\alpha};
\end{eqnarray*}

\item $D^*_p[\widetilde{M}]$ and $D^*_p[\widetilde{T_{R-1}}]$
are Fredholm at $* = p$;

\item\label{v1}
 The $p$-th $L^2$-cohomology of $D^*_p[\widetilde{M}]$ and of
 $D^*_p[\tilde T_{R-1}]$
vanishes. In particular the boundary morphism $\delta^p$ is trivial.
\end{enumerate}
\end{lemma}
\begin{proof}
1.) The isometry of Lemma \ref{isometries}
between $\widetilde{T_{R-1}}$ and
$\widetilde{T_0}$ intertwines
$\Delta_p[\widetilde{T_{R-1}}]$ and
$\Delta_p[\widetilde{T_0}]$.
In particular, spectral
projections which enter in the
definition of the spectral density function are
mapped onto each other. However,
this does not imply that the spectral density
functions are equal because we
have to take regularized dimensions, which
involve the action of the fundamental group,
and the obtained isometries do not
respect the group action.
Note that we explicitly get the regularized
dimension by integrating the trace
of the kernel of the corresponding projection
operator over a fundamental
domain of the diagonal.
The isometry above maps the
fundamental domain of $\widetilde{T_{R-1}}$
onto a subset of the fundamental domain in
$\widetilde{T_0}$ (after suitable
choices).\par\noindent
2.) We conclude from \cite[Proposition 39, Proposition 46]{Lott (1992a)}
that closed hyperbolic manifolds and closed manifolds
with virtually abelian fundamental groups
have positive Novikov-Shubin invariants.
Since these are homotopy invariants and agree with their
combinatorial  counterparts
\cite[Theorem 5.13]{Lott-Lueck (1995)},
\cite{Efremov (1991)},\cite{Gromov-Shubin (1991)},
the same is true for $T_0$. Note that $F(\Delta_p^\perp(\tilde
M),\lambda)=f_{p,\HH^m}(\lambda)\cdot\vol(M)$ depends only on the
homogenous strucure of $\HH^m$ and the volume of the quotient. In
particular, its behavior near zero (given by $f_{p,\HH^m}$ can be
computed using a closed quotient. \par\noindent
3. and 4.)  Because of 2.) $\Delta_p^{\perp}[\widetilde{M}]$ and
$\Delta_p^{\perp}[\widetilde{T_{R-1}}]$ are left-Fredholm.
Moreover, the kernel of $\Delta_p^{\perp}[\widetilde{M}]$
is trivial since $\widetilde{M}$ is $\HH^m$ and for odd $m$
there are no harmonic $L^2$-forms on $\HH^m$ \cite{Dodziuk (1979)}. Also,
$T_{R-1}$ is homotopy equivalent to the torus $T^{m-1}$
which has trivial $L^2$-cohomology.
Since $L^2$-cohomology is a homotopy invariant of
compact manifolds there are no harmonic $L^2$-forms on $\widetilde{T_{R-1}}$.
Now Proposition \ref{Fohne},
Proposition \ref{Fabs},
Proposition \ref{comparision  of Laplace and L} imply 3.)
and the vanishing
of the $L^2$-cohomology of $D^*_p[\widetilde{M}]$  and
$D^*_p[\widetilde{T_{R-1}}]$.
\end{proof}

\begin{corollary}\label{L2_acyclic}
$M_R$ is $L^2$-acyclic for all $R>0$.
\end{corollary}
\begin{proof}
The long weakly exact cohomology sequence of \eqref{ses} together with
Lemma \ref{delta is trivial}.\ref{v1} implies that the $p$-th $L^2$-cohomology
of $D^*_p[\tilde M_R]$ vanishes. Propositions \ref{Fohne}, \ref{Fabs} and
 \ref{comparision of Laplace and L} show
that also the  $L^2$-cohomology as it is usually defined is trivial.
\end{proof}

\begin{lemma} \label{c_E,c_C.c_1 are independent of R}
We can choose the constants
$c_E$, $c_C$ and $c_1$ of Theorem
\ref{short_exact} applied to the sequence
\eqref{ses} at $* = p$ independently of $R$.
\end{lemma}
\begin{proof}
It suffices to  find a constant $C < 0$ independent of $R$
such that for all $R \ge 2$ the following holds
\begin{eqnarray}
\norm{d^p: H^1(\Lambda^p(\widetilde{M_R})) \to
L^2(\Lambda^{p+1}(\widetilde{M_R}))}
& \le & 1; \label{eqn 5.1}
\\
\norm{d^p: H^1(\Lambda^p(\widetilde{E_{R-1}})) \to
L^2(\Lambda^{p+1}(\widetilde{E_{R-1}}))}
& \le & C; \label{eqn 5.2}
\\
\norm{j_{p+1}} & \le & \sqrt{2};
\label{eqn 5.3}
\\
\norm{q_{p+1}} & \le & 1;
\label{eqn 5.4}
\\
\norm{j_{p}} & \le & C;
\label{eqn 5.5}
\\
\norm{q_{p}} & \le & C;
\label{eqn 5.6}
\\
\norm{j_{p+1}^{-1}} & \le & 1;
\label{eqn 5.7}
\\
\norm{j_{p}^{-1}} & \le & 1;
\label{eqn 5.8}
\\
\norm{q_{p}^{-1}} & \le & C;
\label{eqn 5.9}
\end{eqnarray}

We get \eqref{eqn 5.1} directly from the definition of the Sobolev norm
\ref{Sobolev spaces}. We conclude \eqref{eqn 5.2} from the fact that
$E_{R-1}$ is isometrically diffeomorphic to $E_0$ (Lemma \ref{isometries}).
We obtain \eqref{eqn 5.3} from
\[
 \abs{j_{p+1}\omega}^2_{L^2} = \int_{\widetilde{M_R}}\abs{\omega(x)}^2\,dx +
 \int_{\widetilde{E_{R-1}}}\abs{\omega(x)}^2\,dx
 \le \int_{\widetilde{M}}\abs{\omega(x)}^2
+\int_{\widetilde{M}}\abs{\omega(x)}^2 =2\abs{\omega}^2_{L^2}.
\]
and similarly for \eqref{eqn 5.4}. For $j_p$ in \eqref{eqn 5.5} observe
\begin{eqnarray*}
\abs{\omega|_{\widetilde{M_R}}}_{H^1}^2 +
\abs{\omega|_{\widetilde{E_{R-1}}}}_{H^1}^2
& \le  &
 2(\abs{\omega}_{L^2}^2 +\abs{d \omega}^2_{L^2} +
 \abs{\delta\omega}_{L^2}^2)
 + \abs{i^*(*\omega)}^2_{H^{1/2}(\boundary \widetilde{M_R})}
+\abs{i^*(*\omega)}^2_{H^{1/2}(\boundary\widetilde{E_{R-1}})}
\\
& = & 2\abs{\omega}_{H^1}^2 +
\abs{i^*(*\omega)}^2_{H^{1/2}(\boundary \widetilde{M_R})}
+\abs{i^*(*\omega)}^2_{H^{1/2}(\boundary\widetilde{E_{R-1}})}.
\end{eqnarray*}

Therefore, we only have to  deal
with restriction to the boundary.
Choose a cutoff function
$\chi:\R \to [0,1]$ such that for some $0 < \epsilon <1$
$\chi(u) = 0$ for $ \abs{u} \ge  \epsilon$ and
$\chi(u) = 1$ for $ \abs{u} \le \epsilon/2$ holds.
Define a function
$$\chi_R: M \longrightarrow [0,1]$$
which vanishes on $\widetilde{M_0}$ and sends
$(u,x) \in [0,\infty) \times \R^{m-1} = \widetilde{E_0}$ to $\chi(u-R)$.
The support of $\chi_R$ lies in the interior of
$\widetilde{T_{R-1}}\cup \widetilde{T_R}$. Then
\[
 \abs{i^*(*\omega)}^2_{H^{1/2}(\boundary\widetilde{M_R})}
 ~ = ~
\abs{i^*(*\chi_R\omega)}^2_{H^{1/2}(\boundary\widetilde{M_R})}
~ \le ~
C^2_{\widetilde{T_{R-1}}\cup
\widetilde{T_R}}\abs{\chi_R\omega}^2_{H^1(\widetilde{M})}
~ \le ~
C^2_{\widetilde{T_{R-1}}\cup
\widetilde{T_R}} C^2_{\chi_R}\abs{\omega}^2_{H^1}.
\]
The isometries of Lemma \ref{isometries}
which map $\widetilde{T_{R-1}}\cup \widetilde{T_R}$
to $\widetilde{T_0}\cup \widetilde{T_1}$
interchange $\chi_R$ and $\chi_1$. Because the
Sobolev norms  are defined locally in
terms of the geometry, the existence of
these isometries shows that we can choose
$C_{\widetilde{T_{R-1}}\cup \widetilde{T_R}} C_{\chi_R}$
independently of $R$. Since a similiar argument applies to
$\abs{i^*(*\omega)}^2_{H^{1/2}(\boundary\widetilde{E_{R-1}})}$
we get \eqref{eqn 5.5}. The proof of \eqref{eqn 5.6} is similiar.

To show \eqref{eqn 5.7}, for every
$\omega_1\oplus\omega_2\in L^2(\widetilde{M_R})\oplus
L^2(\widetilde{E_{R-1}})$ with
$\omega_1|_{T_{R-1}}=\omega_2|_{T_{R-1}}$  we must find an element
$\omega\in L^2(\widetilde{M})$ with
$j(\omega)=\omega_1\oplus\omega_2$ and
$\abs{\omega}^2_{L^2}\le
\abs{\omega_1}^2_{L^2}+\abs{\omega_2}^2_{L^2}$.
The latter is easy to achieve. Namely, define
\[ \omega(x):=\begin{cases}
\omega_1(x);& x\in \widetilde{M_R}\\
\omega_2(x);& x\in\widetilde{E_{R-1}}
\end{cases} \]
We use the same method to prove \eqref{eqn 5.8}.
It  remains to prove \eqref{eqn 5.9}.

Choose for $R=1$ a bounded operator
\[ \Ex_0: H^1(\Lambda^p(\widetilde{T_0}))
\to H^1(\Lambda^p(\widetilde{E_0})) \]
which satisfies $\Ex_0(\omega)|_{\widetilde{T_0}} = \omega$
for all $\omega \in  H^1(\Lambda^p(\widetilde{T_0}))$.
For arbitrary $R$, define the
corresponding extension operator
\[ \Ex_R: H^1(\Lambda^p(\widetilde{T_R}))
\to H^1(\Lambda^p(\widetilde{E_R})) \]
using $\Ex_0$ and
the isometries of $\widetilde{E_0}$ and $\widetilde{E_{R-1}}$
given in Lemma \ref{isometries}. Since
the $H^1$-norm is defined entirely in
terms of  the Riemannian metric, the norms
of all the extension operators $E_R$ are equal.
We get for $\omega \in  H^1(\Lambda^p(\widetilde{T_0}))$
\begin{eqnarray*}
q_p(0 \oplus \Ex_{R-1}(\omega)) & = & \omega;
\\
\abs{0 \oplus \Ex_{R-1}(\omega)}_{H^1} & \le & 
\norm{\Ex_{R-1}}\abs{\omega}_{H^1}.
\end{eqnarray*}
This implies $\norm{q_p^{-1}} \le \norm{\Ex_{R-1}} = \norm{\Ex_0}$.
This shows \eqref{eqn 5.9} and finishes the proof of
Lemma \ref{c_E,c_C.c_1 are independent of R}.
\end{proof}

\begin{proposition}\label{G_prop}
There is a constant $C$ so that
\[G(\lambda) ~ = ~
F(\Delta^{\perp}_p[\widetilde{T_{0}}],C \cdot \lambda^{1/2}) +
F(\Delta^{\perp}_p[\widetilde{M}],C \cdot \lambda^{1/4}).\]
fulfills the assumptions in Corollary \ref{goal}.
\end{proposition}
\begin{proof}
We conclude from Theorem \ref{short_exact}, equation
\eqref{eqn 5.0}, Lemma \ref{delta is trivial} and Lemma
\ref{c_E,c_C.c_1 are independent of R} that there is a constant
$C > 0$ independent of $R$
such that for all $R \ge 1$ and $0 \le \lambda \le C^{-1}$
\[ F_p(D^*_p[\widetilde{M_R}],\lambda)
~ \le ~
F_p(D^*_p[\widetilde{T_{0}}],C \cdot \lambda^{1/2}) +
F_p(D^*_p[\widetilde{M}],C \cdot \lambda^{1/4}).\]
Now  we apply Proposition \ref{Fohne},
Proposition \ref{Fabs},
Proposition \ref{comparision of Laplace and L}.
One checks easily that all relevant statements also
hold for $\widetilde{M}$ although $M$
itself is not compact because $\widetilde{M}$ is isometric
to $\HH^m$ which is homogeneous. Hence
there is a constant $C > 0$ independent of $R$ such that for all $R \ge 1$
and $0 \le \lambda \le C^{-1}$
\[ F(\Delta_p^{\perp}[\widetilde{M_R}],\lambda)
~ \le ~ G(\lambda).\]
We conclude from
Lemma \ref{delta is trivial}.2
\begin{eqnarray*}
 \int_1^\infty
 \left(\int_0^\epsilon e^{-t\lambda}G(\lambda)\,d\lambda\right)\,dt
& \le &
\int_1^\infty
\left(\int_0^\epsilon e^{-t\lambda}(\lambda^{1/2} + \lambda^{1/4})
\,d\lambda\right)\,dt
\\
& \le &
2\cdot \int_1^\infty
\left(\int_0^\epsilon e^{-t\lambda}\lambda^{1/4}
\,d\lambda\right)\,dt
\\
& \le &
2\cdot
\int_0^\epsilon  \left(\int_1^\infty e^{-t\lambda} \,dt \right)\lambda^{1/4}
\,d\lambda
\\
& \le &
2\cdot
\int_0^\epsilon -e^{-t\lambda}\cdot \lambda^{-3/4}
\,d\lambda
\\
& < & \infty.
\end{eqnarray*}
This proves Proposition \ref{G_prop}.
\end{proof}

This finishes the proof of our main Theorem \ref{Tor_converg} because
of Corollary \ref{goal} and Proposition \ref{G_prop}
for the large $t$ summand and
Lemma \ref{new expression for zetafunction} and Proposition
\ref{dsm} for the small $t$ summand.



\typeout{--------------------- section 7 -----------------------}

\section{$L^2$-analytic torsion and variation of the metric}
\label{L^2-analytic torsion and variation of the metric}

In the next lemma we extend a result of Lott \cite[p.~480]{Lott (1992a)} to
manifolds with boundary.
\begin{lemma}\label{uniform decay}
Let $N$ be a compact manifold, possibly with boundary. Let
$(g_u)_{u\in[0,1]}$ be a continuous family of
Riemannian metrics on $N$. Then
\[ \tr_\Gamma e^{-t\Delta_p^\perp[\tilde N,g_u]}\xrightarrow{t\to\infty} 0
\quad\text{uniformly in $u\in[0,1]$}. \]
Here, $\Delta^\perp$ is $\Delta$ restricted to the orthogonal complement of its
kernel.
\end{lemma}
\begin{proof}
  If $E_\lambda(u)$ is the right continuous spectral family of $\Delta_p^\perp[
\tilde N,g_u]$ and $F(\lambda,u)=\tr_\Gamma E_\lambda(u)$
 we have by \eqref{ew} for
every $\epsilon>0$ and $t\ge 1$
\begin{equation}\label{e1} \tr_\Gamma e^{-t\Delta_p^\perp[\tilde N,g_u]} \le
 \int_0^\epsilon t e^{-t\lambda}F(\lambda,u)\,d\lambda
+ e^{-t\epsilon}F(\epsilon,u)
 + e^{-t\epsilon+\epsilon}\tr_\Gamma e^{-\Delta_p^\perp[\tilde N,u]}.
\end{equation}
By Proposition \ref{comparision of Laplace and L} we have
\begin{equation}\label{e2}
F(\lambda,u) =F_p(L^*_p[\tilde N,g_u],\lambda^2)
+F_{p-1}(L^*_{p-1}[\tilde N,g_u],\lambda^2).
\end{equation}
 The complexes $L^*$ are
defined in Definition \ref{Sobolev cochain complexes}. The identity map
induces a bounded isomorphism between $L^*_p[\tilde N,g_u]$ and $L^*_p[
\tilde N,g_0]$ with norm $a_p(u)$ (the norm is in general different from $1$
because the inner products are different). Denote by $b_p(u)$ the norm of the
inverse. Since $g_u$ is continuous $a:=\sup_{p=0,\dots,m;\;
u\in[0,1]} a_p(u)$ and $b:=
\sup_{p=0,\dots m;\;u\in[0,1]} b_p(u)$ exist.
By Theorem \ref{Gromov-Shubin} we can compare
the spectral density functions in the following way:
\begin{equation}\label{e3}
 F_p(L^*_p[\tilde N,g_u],\lambda)\le F_p(L^*_p[\tilde N,g_0],a^2b^2\lambda).
\end{equation}
Now \eqref{e1}, \eqref{e2} and \eqref{e3} imply
\begin{equation}
  \label{e4}
  \tr_\Gamma e^{-t\Delta_p^\perp[\tilde N,g_u]} \le
  \int_0^\epsilon t e^{-t\lambda}F(ab\lambda,0)\,d\lambda
+ e^{-t\epsilon}F(ab\epsilon,0)
 + e^{-t\epsilon+\epsilon}\tr_\Gamma e^{-\Delta_p^\perp[\tilde N,u]}.
\end{equation}
By its explicit construction, 
the integral kernel of $e^{-\Delta_p^\perp[\tilde N,u]}$ is a continuous
function of $u$ and therefore $\tr_\Gamma e^{-\Delta^\perp[\tilde N,u]}$ is
uniformly bounded in $u$. It follows that $\tr_\Gamma e^{-t\Delta_p^\perp[
\tilde N,g_u]}\xrightarrow{t\to\infty}0$ uniformly in $u$ if we find
$\epsilon>0$ so that $\int_0^\epsilon te^{-t\lambda}F(ab\lambda,0)\,d\lambda
<\infty$. We get for $\epsilon = (ab)^{-1}$, $t\ge 1$
\[\begin{split}
\int_0^\epsilon te^{-t\lambda}F(ab\lambda,0)\,d\lambda &=
 \int_0^{ab\epsilon}\frac{t}{ab}e^{-\lambda t/ab}F(\lambda,0)\,d\lambda\\
& \stackrel{t\ge 1}{\le}
 \int_0^{1}\frac{1}{ab}e^{-\lambda t/ab}F(\lambda,0)\,d\lambda\\
& \le \frac{1}{ab} \cdot \tr_\Gamma e^{-(t/ab)\Delta_p^\perp[\tilde N,g_0]}
<\infty \qed
\end{split}\]
\renewcommand{\qed}{}\end{proof}

\begin{theorem}\label{var_tor1}
Suppose $N^m$ is a compact
manifold, possibly with boundary. Suppose $N$ is of
determinant class. Let $(g_u)_{0\le u\le 1}$
 be a smooth family of Riemannian metrics on $N$.
Let $*_u$ be the Hodge-$*$-operator with respect to the metric $g_u$. Set
$V_p=V_p(g_u):= (\frac{d}{du}*_p(g_u))\circ*_p(g_u)^{-1}$.

The analytic and $L^2$-analytic torsion are smooth functions of $u$, and
\begin{equation}\label{difft}
 \frac{d}{du}|_{u=0}\left(\anTor(N,g_u)-
 \ordanTor(N,g_u)\right) = \sum_p (-1)^{p+1}\left(\tr_\Gamma
V_p|_{\ker \Delta_p(\tilde N)} - \tr V_p|_{\ker\Delta_p(N)}\right).
\end{equation}

\end{theorem}
\begin{proof}
Define
\begin{multline}\label{T_ab}
f(u,\Lambda):=\\
\sum_p (-1)^p p\left(\int_0^\Lambda
 \tr_\Gamma e^{-T\Delta_p(\tilde N)}-\tr e^{-T\Delta_p(N)}
 -\theta(T-1)(b_p^{(2)}-b_p)\frac{dT}{T} -\gamma(b_p^{(2)}-b_p)\right),
\end{multline}
where $\gamma$ is the Euler--Mascheroni constant, $\theta$ the Heaviside step
function and $b_p$ are the ordinary Betti numbers. Lott proves the theorem
 for $\boundary N=\emptyset$
\cite[Proposition 7]{Lott (1992a)}. His proof applies to our situation because
of the following remarks:
\begin{enumerate}
\item
The asymptotic expansions for $\tr_\Gamma e^{-t\Delta_p(\tilde N)}$ and
$\tr e^{-t\Delta_p(N)}$ are equal by \eqref{eqn 1.111} and therefore Lott's
Corollary 5 with proof generalizes to our situation. This implies
  \begin{equation}
    \label{lim}
  f(u,\Lambda)\xrightarrow{\Lambda\to\infty}\anTor(N,g_u)-\ordanTor(N,g_u).
  \end{equation}

\item Set $B(u,T):= \sum_p(-1)^p \left(\tr_\Gamma e^{-T\Delta_p(\tilde N,g_u)}-
\tr e^{-T\Delta_p(N,g_u)}\right)$. Then
\begin{equation}\label{diff}
 \frac{d B}{du} = \sum_p(-1)^{p+1}T\frac{d}{dT}\left(
\tr_\Gamma Ve^{-T\Delta_p(\tilde N)}-\tr Ve^{-T\Delta_p(N)}\right).
\end{equation}
Note that Duhamel's principle works also for the manifold $N$ with boundary,
 so that Lott's proof applies.
The only crucial point is the question wether $\exp{(-T\Delta_p[\tilde
  N,g_u])}$
is differentiable with respect to $u$.
 This follows from the explicit construction of this operator (f.i.~in
 Ray/Singer
\cite[5.4]{Ray-Singer(1971)}) which works also on $\tilde N$ because of
bounded geometry.
\item The principle of not feeling the boundary in
Theorem \ref{kernComp} implies that
$$\frac{d}{dT}\left(\tr_\Gamma Ve^{-T\Delta_p[\tilde N,g_u]}-
\tr V e^{-T\Delta_p[N,g_u]}\right)=
-\tr_\Gamma V\Delta_p[\tilde N]e^{-T\Delta_p[\tilde N]}+\tr V\Delta_p[N]
e^{-T\Delta_p[N]}$$
is bounded for $0\le T\le\Lambda<\infty$ in the same way as it
implies Lemma \ref{asymptotic expansion of the heat kernel} because $V$ is
a local operator. As $g_u$ is
smooth in $u$, the expression is uniformly bounded for $0\le T\le \Lambda$ and
$0\le u\le 1$. Therefore we can integrate \eqref{diff} to get
\begin{eqnarray*}
\frac{d}{du} f(u,\Lambda)
& = & \sum_p(-1)^{p+1}\left(\tr_\Gamma
 V_p e^{-\Lambda\Delta_p[\tilde N]}-\tr V_pe^{-\Lambda\Delta_p[N]}\right).
\end{eqnarray*}
and hence
\begin{eqnarray}
\hspace{-5mm} f(u,\Lambda)
& = & f(0,\Lambda)+ \sum_p(-1)^{p+1}\int_0^u \tr_\Gamma
V_p e^{-\Lambda
\Delta_p[\tilde N,g_v]}-\tr V_pe^{-\Lambda\Delta_p[N,g_v]} \,dv.
\label{diff2}
\end{eqnarray}
\item Since $N$ is of determinant class, $\tr_\Gamma e^{-t\Delta^\perp
[\tilde N,g_u]}
\xrightarrow{t\to\infty} 0$ uniformly
in $u\in[0,1]$ by Lemma \ref{uniform decay}. Here $\Delta^\perp$ is $\Delta$
restricted to the orthogonal complement of its kernel. Because $e^{-t\Delta}
=e^{-t\Delta^\perp}+\pr_{\ker\Delta}$ and $V(g_u)$ is a local operator which is
uniformly bounded in $u$ we can interchange integral and limit in \eqref{diff2}
and conclude with \eqref{lim}
\begin{multline*}
\anTor(N,g_u)-\ordanTor(N,g_u) =\\
= \anTor(N,g_0)-\ordanTor(N,g_0) +
\sum_p(-1)^{p+1}\int_0^u \tr_\Gamma V_p|_{\ker\Delta_p[\tilde N,g_v]} - \tr
V_p|_{\ker\Delta_p[N,g_v]}\,dv .
\end{multline*}
\end{enumerate}
The last equation implies \eqref{difft}.
\end{proof}

\begin{definition}\label{metric anomaly term}
The trace of
 $e^{-t\Delta_p(N,g_0)}V_p$ has an asymptotic expansion for $t\to 0$
because $V_p$ is local. Let $d_p$ be the coefficient of $t^0$ of the boundary
contribution to this asymptotic expansion. This is an integral over a density
on $\boundary N$ which is given locally in terms of the germ of the family of
metrics $g_u$ on $\boundary N$
 (compare Cheeger \cite[p.~278]{Cheeger (1979)}).
\end{definition}

\begin{corollary}\label{metric anomaly}
 In the situation of Theorem \ref{var_tor1} and with Definition
\ref{metric anomaly term}
\[ \left.\frac{d}{du}\right|_{u=0} \anTor(N,g_u) = \sum_p(-1)^{p}\left(d_p -
\tr_\Gamma V_p|_{\ker\Delta_p(\tilde N)}\right) . \]
\end{corollary}
\begin{proof}
  This follows from Theorem \ref{var_tor1} and the computation of
$\frac{d}{du}|_{u=0}\ordanTor(N,g_u)$ by Cheeger \cite[3.27]{Cheeger (1979)}
(note that Cheeger's $\alpha$ is just $-V$).
\end{proof}

\begin{appendix}
  \section{Examples for nontrivial metric anomaly}
  \label{metric_anomaly_example}

In Corollary \ref{metric anomaly}, we extended Cheeger's computation
of the deviation of torsion under variation of the Riemannian metric
on a compact manifold with boundary
from classical analytical torsion to $L^2$-analytic
torsion. Here, we will show that the abstract
correction term (in the acyclic case)
\begin{equation*}
  \sum_p (-1)^p d_p
\end{equation*}
can be nonzero. Note that the formula relating analytic and
topological torsion implies that the correction term is zero as long
as the metric has a product structure near the boundary. Our
examples show that this formula (i.e.\ the extension of the
Cheeger-M\"uller theorem to
manifolds with boundary) is not true in general (Corollary
\ref{counterex}). We use absolute boundary conditions for the
examples (as we have done in the main text), but one can produce
counterexamples with relative boundary conditions exactly in the same way.

Branson/Gilkey \cite{Branson-Gilkey(1990)} explicitely compute  the
first few coefficients of the asymptotic expansion of the heat
operator for manifolds with boundary and local boundary
conditions. We will use these results to give two and three
dimensional examples with nontrivial metric anomaly. Using a product
formula this yields examples also in arbitrary dimensions $>1$.

\subsection{The Examples}
\label{sec:examples}

\subsubsection{Dimension $=2$}
Here, we work with $S^1\times[0,3]$.
On $S^1\times [0,1)$, we choose the following
family of metrics:
\[ g_u= f(x,u)(dx^2+dy^2)\quad\text{with $(y,x)\in S^1\times[0,1)$} . \]
Here $S^1=\R/\Z$ with the induced metric.
We choose $(1,dx,dy,dx\wedge dy)$ as basis for the exterior
algebra. Then
$\abs{dx}=f^{-1/2}=\abs{dy}$ and $(dx,dy)=0$. Consequently
$\abs{dx\wedge dy}=f^{-1}$. For the Hodge-$*$-operator we get the
following
\begin{align*} *1 &= fdx\wedge dy & *dx &= dy & *dy &=- dx
  & *(dx\wedge dy) &= f^{-1}  .
\end{align*}
Consequently (remember that $V_p=\partial */\partial u \cdot *^{-1}$)
\begin{align*}
  V_0 &= -f'f^{-1} &  V_1 &= 0 & V_2 &= f'f^{-1} .
\end{align*}
Here $f'=\partial f/\partial u$.
We extend these metrics to smooth metrics on $S^1\times[0,3]$,which
are constant (in $u$) product metrics on $S^1\times(2,3]$. 
The specific choice does not affect the boundary terms we want to compute.

Now we specify $f(x,u)$: choose $f(x,u)=(1+ux)$.

\begin{proposition}
For the family of metrics $g_u$ on $S^1\times[0,3]$,
 the boundary contribution of
the metric anomaly is nonzero:
\[ d_0-d_1+d_2 \ne 0 . \]
\end{proposition}
\begin{proof}
  We have to evaluate the corresponding expressions in
  \cite[7.2]{Branson-Gilkey(1990)}. $d_p$ is the coefficient of $t^0$,
  which in dimension two is the second nontrivial coefficient.
Therefore
\[ d_p = (24\pi)^{-1}\left( \int_N (6 V_p E+F_p\tau) +
  \int_{\boundary N}(2V_p L_{aa}+3\psi_p (V_p)_{;N} +12 V_p S)
\right) . \]
 The integral over $N$ does not matter, because we already know that
 there is no interior contribution to the metric anomaly (in the
 acyclic case), so the alternating sum of these summands has to
 vanish and we do not have to specify the terms in the integrant.

 Since
  $f'f^{-1}=0$ at the boundary, $V_p=0$ at the boundary, and the
  boundary contribution in $d_p$ reduces to
\[ (8\pi)^{-1} \int_{\boundary N}(\psi_p(V_p)_{;N}) . \]
Here $(V_p)_{;N}$
denotes covariant differentiation in normal direction to the
boundary. Since $\partial/\partial x$ is the inward pointing unit
normal to the boundary: $\phi_{;N}=\partial\phi/\partial x$. In
particular
\begin{align*}
  (V_0)_{;N}&=-(1+ux)^{-1} + ux(1+ux)^{-2}=-(V_2)_{;N} &\text{and } (V_1)_{;N}& =0 .
\end{align*}
 $\psi_p$
is the dimension of the subspace which fulfills Neumann boundary
conditions minus the dimension of the subspace with Dirchlet boundary
conditions, i.e.
$\psi_0=1$, where Neumann boundary conditions are in
effect, and $\psi_2=-1$, because on $\Lambda^2$ we have to impose Dirichlet
boundary conditions.

Taking everything together yields (since $x=0$ at the boundary)
\[ d_0-d_1+d_2 = -(4\pi)^{-1}\int_{S^1} 1\; dy
=-(4\pi)^{-1}\ne 0 . \qed\]
\renewcommand{\qed}{}
\end{proof}

\subsubsection{Dimension $=3$}
To produce a three dimensional example, we make a similar Ansatz on
$[0,3]\times T^2$:
\[ g_u = f(u,x)^2 (dx^2 +dy^2+dz^2) \qquad \text{with
  }f(u,x)=(1+x+ux) . \]
Here $(x,y,z)\in [0,1)\times\R^2/\Z^2$, and we give $T^2=\R^2/\Z^2$ the
quotient metric.
The boundary is given by $x=0$. Again we extend these metrics, so that
$(2,3]\times T^2$ gets the standard product metric.

Computations as above give
\begin{align*}
  *1&=f^3 dx\wedge dy\wedge dz\\
  *dx&= f\,dy\wedge dz & *dy &= f\,dz\wedge dx & *dz &= f\,dx\wedge
  dy\\
  *(dy\wedge dz) &=f^{-1}dx & *(dz\wedge dx) &=f^{-1}dy & *(dx\wedge dy)
  &= f^{-1}dz\\
  *(dx\wedge dy\wedge dz)&=f^{-3} .
\end{align*}
It follows
\begin{align}\label{Vi}
V_0 &=   -3 f'f^{-1}  & V_1 &= -f'f^{-1}\\
V_2 &= f'f^{-1} & V_3 &= 3f'f^{-1}\qquad\text{with
  }f'f^{-1}=x(1+x+ux)^{-1}\notag .
\end{align}
In dimension three, the coefficient of $t^0$ is the third nontrivial
coefficient. Branson/Gilkey \cite[7.2]{Branson-Gilkey(1990)} compute
this
coefficient as follows:
\begin{equation*}
  \begin{split}
    1536\pi\cdot d_p = \int_{\boundary N} & V_p\times T\\
                & +(V_p)_{;N}((6\psi^p_N+30\psi^p_D)k+96 S_p)\\
                & + 24\psi_p (V_p)_{;;N} . 
  \end{split}
\end{equation*}
 Here $T$ is a complicated expression in terms of the geometry but
 $V_p(0)=0$ due to our choice of $f$. Since $\partial_x$ is the unit
 normal vector to the boundary, as above
\begin{equation}\label{coeff}
 d_p = (256\pi)^{-1}\int_{\boundary N} 4\psi_p 
\underbrace{(V_p)_{;\partial_x\partial_x}}_{\text{iterated covariant
    derivative in normal direction}} +
\frac{\partial V_p}{\partial x} \left( (\psi^p_N+5\psi^p_D)k 
+16 S_p\right)\;dy\,dz  . 
\end{equation}
Here $k$ is the mean curvature of the boundary (the trace of the
second fundamental form). $\psi^p_N$ is the trace of the projection
onto the subspace of $\Lambda^p$, where
Neumann boundary conditions are in effect, $\psi^p_D$ is the trace 
of the projection onto the subspace with Dirichlet boundary conditions and
$\psi=\psi_N-\psi_D$.
Explicitely we get: 
\begin{tabular}[t]{l|llll}
$p$ & 0 & 1 & 2 & 3\\\hline
$\psi_N^p$ & 1 & 2 & 1 & 0\\
$\psi_D^p$ & 0 & 1 & 2 & 1\\
$\psi_p$  & 1 & 1 & -1 & -1 .
\end{tabular}

Together with \eqref{Vi} this easily implies  that in our example most
of the summands cancel. We are left with
\[ d_0-d_1+d_2-d_3 = (16\pi)^{-1}\sum_{p=0}^3(-1)^p\int_{\boundary N}
\frac{\partial V_p}{\partial x} S_p .
\]
It remains to identify $S_p$. This
 is the trace of the (local) operator $A:=\partial/\partial x -
\nabla_N$ restricted to the subspace of $\Lambda^p$ where we impose
Neumann boundary condition, i.e.~restriction
to all of $\Lambda^0$, to the span of
$dy$ and $dz$ in $\Lambda^1$,  to the span of $dy\wedge dz$ in
$\Lambda^2$, and to  zero in $\Lambda^3$. 
Therefore $S_0=0=S_3$. For $3$-manifolds, $A_1$ and $A_2$ are given in
the proof of \cite[7.2]{Branson-Gilkey(1990)} (table on the bottom of
p.~267). We get: $S_1=-k=S_2$.
In our situation, using the fact the $\partial_x$ is normal to the
boundary and has unit length there, that the Levi-Civita connection is
compatible with the
metric and torsion free
\[\begin{split}
k &=(\nabla_y\partial_y,\partial_x)+(\nabla_z\partial_z,\partial_x)
=(\partial_y,\nabla_y\partial_x)+(\partial_z,\nabla_z\partial_x)\\
&=(\partial_y,\nabla_x\partial_y)+(\partial_z,\nabla_x\partial_z) =
\nabla_x(\partial_y,\partial_y)/2+\nabla_x(\partial_z,\partial_z)/2\\
&=\partial f^2/\partial x = 2(1+x+ux)(1+u) .
\end{split} \]
It follows with \eqref{Vi} (remember $x=0$ on $\boundary N$)
\[  d_0-d_1+d_2-d_3 = -(8\pi)^{-1}\int_{T^2} \frac{\partial
f'f^{-1}}{\partial x} k
\,dy\,dz = -(4\pi)^{-1}(1+u) \ne 0 .
\]
Here we used
\[  \frac{\partial
f'f^{-1}}{\partial x} = (1+x+ux)^{-1} - x(1+u)(1+x+ux)^{-2} 
\stackrel{x=0}{=} 1 . \]

\subsubsection{Arbitrary dimensions}
\begin{proposition}\label{metric_ano_ex}
For every dimension $m>1$ we find a compact manifold $N$ with boundary with a
family of metrics $g_u$ on $N$ so that the boundary contribution to
the metric anomaly is nonzero.
\end{proposition}
\begin{proof}
  For $m=2$ and $m=3$, we have shown that we can use $N=S^1\times [0,3]$
  or $N=S^1\times S^1\times[0,3]$, respectively. For higher
  dimensions, use $N\times S^{2k}$ with suitable $k\in\N$ with the
  corresponding family of product metrics. We use the product formula
  \[ \ordanTor(N\times S^{2k}) = \ordanTor(N)\chi(S^{2k}) +
  \chi(N)\ordanTor(S^{2k}) = 2\ordanTor(N) \]
Using an acyclic representation on $N$ we get nontrivial metric
anomaly of $N\times S^{2k}$ which is entirely due to the boundary. For
$L^2$-analytic torsion observe that $N$ is $L^2$-acyclic so that the
same argument applies.
\end{proof}

\begin{corollary}\label{counterex}
  For every dimension $m>1$ we find compact Riemannian manifolds $M^m$
  with boundary
  which are $L^2$-acyclic, so that the difference between
  $L^2$-analytic and $L^2$-topological torsion is different from
  $\chi(\boundary M)(\ln 2)/2$. Similarly, we find acyclic finite
  dimensional orthogonal representations of $\pi_1(M)$ so that the
  same statement holds for the corresponding classical analytic and
  topological torsion.

  In other words, we get counterexamples for the extension of the
  Cheeger-M\"uller theorem (and its $L^2$-counterpart) to arbitrary compact
  Riemannian manifolds with boundary.
\end{corollary}
\begin{proof}
  The manifolds of Proposition \ref{metric_ano_ex} do the job
  for metrics $g_{u_0}$ with $u_0\ne 0$ sufficiently
  small. This follows from the validity of the Cheeger-M\"uller
  theorem  for the product metric case $g_0$ (by
  \cite{Lueck(1993),Burghelea-Friedlander-Kappeler (1996a)}), and the
  anomaly formulas
  \cite[3.27]{Cheeger (1979)} or
  Corollary \ref{metric anomaly}, respectively.
\end{proof}
\end{appendix}


\typeout{--------------------- References -----------------------}

\begin{center}
 Current address\\
Wolfgang L\"uck  and Thomas Schick\\
Fachbereich Mathematik und Informatik\\
Westf\"alische Wilhelms-Universit\"at M\"unster\\
Einsteinstr. 62\\
48149 M\"unster\\
Bundesrepublik Deutschland\\
email: lueck/schickt@math.uni-muenster.de\\
FAX: 0251 8338370\\
internet: http://wwwmath.uni-muenster.de/math/u/lueck/
\end{center}

\begin{center} Version of \today \end{center}

\end{document}